\documentclass[11pt, DIV10,a4paper]{article}

\usepackage{natbib}
\newcommand {\ctn}{\citet} 
\usepackage{graphicx,subfigure,amsmath,latexsym,amssymb}
\usepackage{float,epsfig,multirow,rotating,times}
\usepackage{upgreek,wrapfig}
\usepackage{comment}

\usepackage{algorithm}

\newcommand{\btheta}{\boldsymbol{\theta}}

\newcommand{\bbeta}{\boldsymbol{\beta}}

\newcommand{\bxi}{\boldsymbol{\xi}}
\newcommand{\bGamma}{\boldsymbol{\Gamma}}
\newcommand{\bTheta}{\boldsymbol{\Theta}}

\newcommand{\bSigma}{\boldsymbol{\Sigma}}
\newcommand{\bPsi}{\boldsymbol{\Psi}}
\newcommand{\bpsi}{\boldsymbol{\psi}}

\newcommand{\bmu}{\boldsymbol{\mu}}

\newcommand{\bvarphi}{\boldsymbol{\varphi}}

\newcommand{\bOmega}{\boldsymbol{\Omega}}

\newcommand{\bB}{\boldsymbol{B}}
\newcommand{\bC}{\boldsymbol{C}}

\newcommand{\bA}{\boldsymbol{A}}
\newcommand{\bS}{\boldsymbol{S}}
\newcommand{\bt}{\boldsymbol{t}}

\newcommand{\bR}{\boldsymbol{R}}

\newcommand{\bZ}{\boldsymbol{Z}}
\newcommand{\bz}{\boldsymbol{z}}

\newcommand{\bzero}{\boldsymbol{0}}

\newtheorem{theorem}{Theorem}

\newtheorem{lemma}[theorem]{Lemma}

\newcommand{\boldcal}[1]{\mbox{\boldmath{$\mathcal #1 $}}}

\newcommand{\bM}{\mathbf M}
\newcommand{\bN}{\mathbf N}

\numberwithin{equation}{section}
\numberwithin{algo}{section}
\numberwithin{table}{section}
\numberwithin{figure}{section}

\usepackage{setspace}
\usepackage[mathscr]{euscript}
\usepackage[margin=1in]{geometry}
\begin{document}

\normalsize

\title{\vspace{-0.8in}
On Asymptotic Inference in Stochastic Differential Equations with Time-Varying Covariates}
\author{Trisha Maitra and Sourabh Bhattacharya\thanks{
Trisha Maitra is a PhD student and Sourabh Bhattacharya 
is an Associate Professor in
 Interdisciplinary Statistical Research Unit, Indian Statistical
Institute, 203, B. T. Road, Kolkata 700108.
Corresponding e-mail: sourabh@isical.ac.in.}}
\date{\vspace{-0.5in}}
\maketitle%

\begin{abstract}
In this article, we introduce a system of stochastic differential equations ($SDE$s) consisting of time-dependent covariates and consider both fixed and random effects set-ups. We also
allow the functional part associated with the drift function to depend upon unknown parameters. In this general set-up of SDE system we establish consistency and asymptotic normality of the
M LE through verification of the regularity conditions required by existing relevant theorems. Besides, we consider the Bayesian approach to learning about the population parameters, and prove
consistency and asymptotic normality of the corresponding posterior distribution. We supplement our theoretical investigation with simulated and real data analyses, obtaining encouraging results
in each case.
\\[2mm]
{\it {\bf Keywords:} 
Asymptotic normality; Functional data; Gibbs sampling; Maximum likelihood estimator; Posterior consistency; Random effects.
}
 
\end{abstract}

\section{Introduction}
\label{sec:intro}

Systems of stochastic differential equations ($SDE$s) are appropriate for modeling situations 
where ``within" subject variability is caused by some random component
varying continuously in time; hence, $SDE$ systems are also appropriate for modeling functional data (see, for example, \ctn{Zhu11}, \ctn{Ramsay05} for some connections between
$SDE$ and functional data analysis). When suitable time-varying covariates are available, it is
then appropriate to incorporate such information in the $SDE$ system.
Some examples of statistical applications of $SDE$-based models with time-dependent covariates 
are \ctn{Zita11}, \ctn{Overgaard05}, \ctn{Leander15}. 

However, systems of $SDE$ based models consisting of time-varying covariates seem to be
rare in the statistical literature, in spite of their importance, and their asymptotic properties are hitherto unexplored.
Indeed, although asymptotic inference in single, fixed effects $SDE$ models without covariates has been considered in the literature as time tends to infinity (see, for example, \ctn{Bishwal08}),
asymptotic theory in systems of $SDE$ models is rare, and 
so far only random effects $SDE$ systems 
without covariates have been considered, as $n$, the number of subjects (equivalently, the number of $SDE$s in the system), tends to infinity 
(\ctn{Maud12}, \ctn{Maitra14a}, \ctn{Maitra14b}). Such models are of the following form:
\begin{equation}
d X_i(t)=b(X_i(t),\phi_i)dt+\sigma(X_i(t))dW_i(t),\quad\mbox{with}\quad X_i(0)=x^i,~i=1,\ldots,n,
\label{eq:sde_basic1}
\end{equation}
where, for $i=1,\ldots,n$, $X_i(0)=x^i$ is the initial value of the stochastic process $X_i(t)$, which
is assumed to be continuously observed on the time interval $[0,T_i]$; $T_i>0$ assumed to be known.
The function $b(x,\phi)$, which is the drift function, is a known, real-valued function on $\mathbb R\times\mathbb R^d$
($\mathbb R$ is the real line and $d$ is the dimension), and the function $\sigma:\mathbb R\mapsto\mathbb R$ is the known
diffusion coefficient.
The $SDE$s given by (\ref{eq:sde_basic1}) are driven by independent standard Wiener processes $\{W_i(\cdot);~i=1,\ldots,n\}$,
and $\{\phi_i;~i=1,\ldots,n\}$, which are to be interpreted as the random effect parameters associated
with the $n$ individuals, which are assumed by \ctn{Maud12} to be independent of the Brownian motions and
independently and identically distributed ($iid$) random variables with some common distribution.

For the sake of convenience
\ctn{Maud12} (see also \ctn{Maitra14a} and \ctn{Maitra14b}) assume $b(x,\phi_i)=\phi_ib(x)$.
Thus, the random effect is a multiplicative factor of the drift function; also, the function $b(x)$ is assumed
to be independent of parameters.
In this article, we generalize the multiplicative factor to include time-dependent covariates; we also allow
$b(x)$ to depend upon unknown parameters.

Notably, such model extension has already been provided in \ctn{Maitra15} and \ctn{Maitra15b}, but their
goal was to develop asymptotic theory of Bayes factors for comparing systems of $SDE$s, with or without
time-varying covariates, emphasizing, when time-varying covariates are present, simultaneous asymptotic selection of covariates 
and part of the drift function free of covariates, using Bayes factors. 

In this work, we deal with parametric asymptotic inference, both frequentist and Bayesian, in the context of our extended
system of $SDE$s. We consider, separately, fixed effects as well as random effects. The fixed effects
set-up ensues when coefficients associated with the covariates are the same for all the subjects. On the other hand,
in the random effects set-up, the subject-wise coefficients are assumed to be a random sample from some distribution
with unknown parameters. 

It is also important to distinguish between the $iid$ situation and the 
independent but non-identical case (we refer to the latter as non-$iid$) that we consider.
The $iid$ set-up is concerned with the case where the initial values $x^i$ and time limit $T_i$ 
are the same for all $i$, and the coefficients associated with the covariates are zero, that is, there
are no covariates suitable for the $SDE$-based system.
This set-up, however, does not reduce to the $iid$ set-up considered in \ctn{Maud12}, \ctn{Maitra14a} and \ctn{Maitra14b}
because in the latter works $b(x)$ was assumed to be free of parameters, while in this work we allow
this function to be dependent on unknown parameters. The non-$iid$ set-up assumes either or both of the
following: presence of appropriate covariates and that $x^i$ and $T_i$ are
not the same for all the subjects.

In the classical paradigm, we investigate consistency and asymptotic normality of the maximum likelihood estimator ($MLE$) of the
unknown parameters which we denote by $\btheta$, and in the Bayesian framework we study consistency and asymptotic 
normality of the Bayesian posterior distribution
of $\btheta$. In other words, we consider prior distributions $\pi(\btheta)$ of $\btheta$ and study the
properties of the corresponding posterior
\begin{equation}
\pi_n(\btheta|X_1,\ldots,X_n)=\frac{\pi(\btheta)\prod_{i=1}^n f_i(X_i|\btheta)}
{\int_{\bpsi\in\bTheta}\pi(\bpsi)\prod_{i=1}^n f_i(X_i|\bpsi)d\bpsi}
\label{eq:posterior1}
\end{equation}
as the sample size $n$ tends to infinity. Here $f_i(\cdot|\btheta)$ is the density corresponding to the $i$-th individual and 
$\bTheta$ is the parameter space. 

In what follows, after introducing our model and the associated likelihood in Section \ref{sec:model}, 
we investigate asymptotic properties of $MLE$ in the $iid$  and non-$iid$ contexts in 
Sections \ref{sec:consistency_iid} and \ref{sec:consistency_non_iid} respectively. Then, in 
Sections \ref{sec:b_consistency_iid} and \ref{sec:b_consistency_non_iid} we 
investigate asymptotic properties of the posterior in the $iid$ and non-$iid$ cases, respectively. 
In Section \ref{sec:another} we consider the random effects set-up and provide necessary discussion to point towards
validity of the corresponding asymptotic results.  
We demonstrate the applicability of our developments to practical and finite-sample contexts using simulated and real data analyses
in Sections \ref{sec:simulated_data} and \ref{sec:truedata}, respectively.
We summarize our contribution and provide further discussion in Section \ref{sec:conclusion}.


Notationally, ``$\stackrel{a.s.}{\rightarrow}$", ``$\stackrel{P}{\rightarrow}$" and ``$\stackrel{\mathcal L}{\rightarrow}$"
denote convergence ``almost surely", ``in probability" and ``in distribution", respectively.

\section{The $SDE$ set-up}
\label{sec:model}
We consider the following system of $SDE$ models for $i=1,2,\ldots,n$:
\begin{equation}
d X_i(t)=\phi_{i,\bxi}(t)b_{\bbeta}(X_i(t))dt+\sigma(X_i(t))dW_i(t)
\label{eq:sde1}
\end{equation}
where $X_i(0)=x^i$ is the initial value of the stochastic process $X_i(t)$, which is assumed to be continuously 
observed on the time interval $[0,T_i]$; $T_i>0$ for all $i$ and assumed to be known. In the above, $\phi_{i,\bxi}$ is the parametric function consisting of
the covariates and the unknown coefficients $\bxi$ associated with the covariates, and $b_{\bbeta}$ is a parametric function known up to the parameters $\bbeta$.

\subsection{Incorporation of time-varying covariates}
\label{subsec:covariates}
We assume that $\phi_{i,\bxi}(t)$ 
has the following form:  
\begin{equation}
\phi_{i,\bxi}(t)=\phi_{i,\bxi}(\bz_i(t))
=\xi_{0}+\xi_{1}g_1(z_{i1}(t))+\xi_{2}g_2(z_{i2}(t))+\cdots+\xi_{p}g_p(z_{ip}(t)),
\label{eq:phi_model}
\end{equation}
where $\bxi=(\xi_{0},\xi_{1},\ldots,\xi_{p})$ is a set of real constants 
 and $\bz_i(t)=(z_{i1}(t),z_{i2}(t),\ldots,z_{ip}(t))$ is the set of available covariate information 
corresponding to the $i$-th individual, depending upon time $t$. We assume that $\bz_i(t)$ is continuous in $t$,
$z_{il}(t)\in \bZ_l$  where $\bZ_l$ is compact and 
$g_l:\bZ_l\rightarrow \mathbb R$ is  continuous, for $l=1,\ldots,p$. Let $\boldcal Z=\bZ_1\times\cdots\times\bZ_p$.
We let $$\mathfrak Z=\left\{\bz(t)\in\boldcal Z:t\in[0,\infty)~\mbox{such that}~\bz(t)~\mbox{is continuous in}~t\right\}.$$
Hence, $\bz_i\in\mathfrak Z$ for all $i$. The function $b_{\bbeta}$ is multiplicative part of the drift function free of the covariates. Note that $\bxi$ consists of $p+1$ parameters. 
Assuming that $\bbeta\in\mathbb R^q$, where $q\geq 1$, it follows that our parameter set 
$\btheta=(\bbeta,\bxi)$ belongs to the $(p+q+1)$-dimensional real space $\mathbb R^{p+q+1}$. 
The true parameter set is denoted by $\btheta_0$.

\subsection{Likelihood}
\label{subsec:likelihood}
We first define the following quantities:
\begin{equation}
U_{i,\btheta} =\int_0^{T_i}\frac{\phi_{i,\bxi}(s)b_{\bbeta}(X_i(s))}{\sigma^2(X_i(s))}dX_i(s),
\quad\quad V_{i,\btheta} =\int_0^{T_i}\frac{\phi^2_{i,\bxi}(s)b_{\bbeta}^2(X_i(s))}{\sigma^2(X_i(s))}ds
\label{eq:u_v}
\end{equation}
for $i=1,\ldots,n$. 

Let $\bC_{T_i}$ denote the space of real continuous functions $(x(t), t \in [0,T_i ])$ defined on $[0,T_i ]$, 
endowed with the $\sigma$-field $\mathcal C_{T_i}$ associated with the topology of uniform convergence 
on $[0,T_i ]$. We consider the distribution $P^{x_i,T_i,\bz_i}$ on $(C_{T_i} ,\mathcal C_{T_i})$ 
of $(X_i (t), t\in [0,T_i ])$ given by (\ref{eq:sde1}) 
We choose the dominating measure $P_i $ as the distribution of (\ref{eq:sde1}) 
with null
drift. So, 
\begin{align}
\frac{dP^{x_i,T_i,\bz_i}}{dP_i}=f_{i}(X_i|\btheta)
&=\exp\left(U_{i,\btheta}-\frac{V_{i,\btheta}}{2}\right).
\label{eq:densities}
\end{align}

\section{Consistency and asymptotic normality of $MLE$ in the $iid$ set-up}
\label{sec:consistency_iid}

In the $iid$ set up we have $x_i=x$ and $T_i=T$ for all $i=1,\ldots,n$. Moreover, the covariates 
are absent, that is, $\xi_{i}=0$ for $i=1,\ldots,p.$ 
Hence, the resulting parameter set in this case is $\btheta=(\bbeta,\xi_{0}).$ 

\subsection{Strong consistency of $MLE$}
\label{subsec:MLE_consistency_iid}

Consistency of the $MLE$ under the $iid$ set-up can be verified by validating the regularity
conditions of the following theorem (Theorems 7.49 and 7.54 of \ctn{Schervish95}); for our purpose
we present the version for compact parameter space.
\begin{theorem}[\ctn{Schervish95}]
\label{theorem:theorem1}
Let $\{X_n\}_{n=1}^{\infty}$ be conditionally $iid$ given $\btheta$ with density $f_1(x|\btheta)$
with respect to a measure $\nu$ on a space $\left(\mathcal X^1,\mathcal B^1\right)$. Fix $\btheta_0\in\bTheta$, and define,
for each $\bM\subseteq\bTheta$ and $x\in\mathcal X^1$,
\[
Z(\bM,x)=\inf_{\bpsi\in \bM}\log\frac{f_1(x|\btheta_0)}{f_1(x|\bpsi)}.
\]
Assume that for each $\btheta\neq\btheta_0$, there is an open set $\bN_{\btheta}$ such that $\btheta\in \bN_{\btheta}$ and
that $E_{\btheta_0}Z(\bN_{\btheta},X_i)> -\infty$. 
Also assume that $f_1(x|\cdot)$ is continuous at $\btheta$ 
for every $\btheta$, a.s. $[P_{\btheta_0}]$. Then, if $\hat\btheta_n$ is the $MLE$ of 
$\btheta$ corresponding to $n$ observations, 
it holds that $\underset{n\rightarrow\infty}{\lim}~\hat\btheta_n=\btheta_0$, a.s. $[P_{\btheta_0}]$.
\end{theorem}

\subsubsection{Assumptions}
\label{subsubsec:assumptions_consistency_iid}

We assume the following conditions:
\begin{itemize}
\item[(H1)] The parameter space $\bTheta=\mathfrak B\times\bGamma$ such that $\mathfrak B$ 
and $\bGamma$ are compact.
\item[(H2)] 
$b_{\bbeta}(\cdot)$ and  $\sigma(\cdot)$ are $C^1$ (differentiable with continuous first derivative) 
on $\mathbb R$ and satisfy
$b^2_{\bbeta}(x)\leq K_1(1+x^2+\|\bbeta\|^2)$  
and $\sigma^2(x)\leq K_2(1+x^2)$
for all $x\in\mathbb R$, for some $K_1, K_2>0$. Now, due to (H1) 
the latter boils down to assuming
$b^2_{\bbeta}(x)\leq K(1+x^2)$  
and $\sigma^2(x)\leq K(1+x^2)$
for all $x\in\mathbb R$, for some $K>0$. 
\end{itemize}
We further assume:
\begin{itemize}
\item[(H3)] For every $x$, let $b_{\bbeta}$ be continuous in $\bbeta=(\beta_{1},\ldots,\beta_{q})$ 
and moreover, for $j=1,\ldots,q$, 
$$\underset{\bbeta\in\mathfrak B}{\sup}~\frac{\left|\frac{\partial b_{\bbeta}(x)}{\partial\beta_{j}}\right|}{\sigma^2(x)}
\leq c\left(1+|x|^{\gamma}\right),$$ 
for some $c>0$ and $\gamma\geq 0$.
\item[(H4)] 
\begin{equation}
\frac{b^2_{\bbeta}(x)}{\sigma^2(x)}
\leq K_{\bbeta} \left(1+x^2+\|\bbeta\|^2\right),
\label{eq:H4_prime_1}
\end{equation}
where $K_{\bbeta}$ is continuous in $\bbeta$.
\end{itemize}

\subsubsection{Verification of strong consistency of $MLE$ in our $SDE$ set-up}
\label{subsubsec:MLE_consistency_iid}

To verify the conditions of Theorem \ref{theorem:theorem1} in our case, note that assumptions 
(H1) -- (H4)  clearly imply continuity of the density $f_1(x|\btheta)$ in the same way as 
the proof of Proposition 2 of \ctn{Maud12}. It follows that $U_{\btheta}$ and $V_{\btheta}$ are continuous in
$\btheta$, the property that we use in our proceedings below.

Now consider,
\begin{align}
Z(\bN_{\btheta},X)=&\inf_{\btheta_1\in \bN_{\btheta}}\log\frac{f_1(X|\btheta_0)}{f_1(X|\btheta_1)}\notag\\
=& \left(U_{\btheta_0}-\frac{V_{\btheta_0}}{2}\right)-
\inf_{\btheta_1\in \bN_{\btheta}}\left(U_{\btheta_1}-\frac{V_{\btheta_1}}{2}\right)\notag\\
\geq& \left(U_{\btheta_0}-\frac{V_{\btheta_0}}{2}\right)-
\inf_{\btheta_1\in \bar\bN_{\btheta}}\left(U_{\btheta_1}-\frac{V_{\btheta_1}}{2}\right)\notag\\
=& \left(U_{\btheta_0}-\frac{V_{\btheta_0}}{2}\right)-\left(U_{\btheta_1^*(X)}-\frac{V_{\btheta_1^*(X)}}{2}\right),
\label{eq:ratio1}
\end{align}
where $\bN_{\btheta}$ 
is an appropriate open subset of the relevant compact parameter space, 
and $\bar\bN_{\btheta}$ is a closed
subset of $\bN_{\btheta}$. 
The infimum of $\left( U_{\btheta_1}-\frac{V_{\btheta_1}}{2}\right)$ is attained at 
$\btheta^*_1=\btheta^*_1(X)\in \bar\bN_{\btheta}$ due to continuity of $U_{\btheta}$ and $V_{\btheta}$
in $\btheta$.

Let $E_{\btheta_0}(V_{\btheta_1})=\breve V_{\btheta_1}$ and $E_{\btheta_0}(U_{\btheta_1})=\breve U_{\btheta_1}$. 
From Theorem 5 of \ctn{Maitra14a} it follows that the above expectations are continuous in $\btheta_1$.
Using this we obtain
\begin{align}
E_{\btheta_0}\left(U_{\btheta_1^*(X)}-\frac{V_{\btheta_1^*(X)}}{2}\right)
&=E_{\btheta^*_1(X)|\btheta_0}E_{X|\btheta^*_1(X)=\bvarphi_1,\btheta_0}
\left(U_{\btheta_1^*(X)=\bvarphi_1}-\frac{V_{\btheta_1^*(X)=\bvarphi_1}}{2}\right)\notag\\
&=E_{\btheta^*_1(X)|\btheta_0}\left(\breve U_{\bvarphi_1}-\frac{\breve V_{\bvarphi_1}}{2}\right)\notag\\
&\leq E_{\btheta^*_1(X)|\btheta_0}\left[\sup_{\bvarphi_1\in \bar\bN_{\btheta}}\left(\breve U_{\bvarphi_1}-\frac{\breve V_{\bvarphi_1}}{2}\right)\right]\notag\\
&=E_{\btheta^*_1(X)|\btheta_0}\left(\breve U_{\bvarphi_1^*}-\frac{\breve V_{\bvarphi_1^*}}{2}\right)\notag\\
&=\left(\breve U_{\bvarphi_1^*}-\frac{\breve V_{\bvarphi_1^*}}{2}\right),
\label{eq:uvexp_new1}
\end{align}
where $\bvarphi^*_1\in \bar\bN_{\btheta}$ is where the supremum of 
$\left(\breve U_{\bvarphi_1}-\frac{\breve V_{\bvarphi_1}}{2}\right)$ is achieved. 
Since $\bvarphi^*_1$ is independent of $X$,  the last step (\ref{eq:uvexp_new1}) follows.

Noting that $E_{\btheta_0}\left(U_{\btheta_0}-\frac{V_{\btheta_0}}{2}\right)$ and 
$\left(\breve U_{\bvarphi_1^*}-\frac{\breve V_{\bvarphi_1^*}}{2}\right)$ are finite due to 
Lemma 1 of \ctn{Maitra15}, it follows that $E_{\btheta_0}Z(\bN_{\btheta},X)> -\infty$.
Hence, $\hat\btheta_n\stackrel{a.s.}{\rightarrow}\btheta_0$ $[P_{\btheta_0}]$, as $n\rightarrow\infty$. 
We summarize the result in the form of the following theorem:
\begin{theorem}
\label{theorem:consistency_iid}
Assume the $iid$ setup and conditions (H1) -- (H4). 
Then the $MLE$ is strongly consistent
in the sense that as $n\rightarrow\infty$,
$\hat\btheta_n\stackrel{a.s.}{\rightarrow}\btheta_0$~ $[P_{\btheta_0}]$.
\end{theorem}

\subsection{Asymptotic normality of $MLE$}
\label{subsec:MLE_iid_normality}

To verify asymptotic normality of $MLE$ we invoke the following theorem provided in \ctn{Schervish95} (Theorem 7.63):
\begin{theorem}[\ctn{Schervish95}]
\label{theorem:theorem2}
Let $\bTheta$ be a subset of $\mathbb R^{p+q+1}$, and let $\{X_n\}_{n=1}^{\infty}$ be conditionally $iid$ given $\btheta$
each with density $f_1(\cdot|\btheta)$. Let $\hat\btheta_n$ be an $MLE$. Assume that
$\hat\btheta_n\stackrel{P}{\rightarrow}\btheta$ under $P_{\btheta}$ for all $\btheta$. Assume that $f_1(x|\btheta)$
has continuous second partial derivatives with respect to $\btheta$ and that differentiation can be passed under the
integral sign. Assume that there exists $H_r(x,\btheta)$ such that, for each $\btheta_0\in int(\bTheta)$ and each
$k,j$,
\begin{align}
\sup_{\|\btheta-\btheta_0\|\leq r}\left\vert\frac{\partial^2}{\partial\theta_k\partial\theta_j}\log f_1(x|\btheta_0)
-\frac{\partial^2}{\partial\theta_k\partial\theta_j}\log f_1(x|\btheta)\right\vert\leq H_r(x,\btheta_0),
\label{eq:h1}
\end{align}
with
\begin{equation}
\lim_{r\rightarrow 0}E_{\btheta_0}H_r\left(X,\btheta_0\right)=0.
\label{eq:h2}
\end{equation}
Assume that the Fisher information matrix $\mathcal I(\btheta)$ is finite and non-singular. Then, under $P_{\btheta_0}$,
\begin{equation}
\sqrt{n}\left(\hat\btheta_n-\btheta_0\right)\stackrel{\mathcal L}{\rightarrow}N\left(\bzero,\mathcal I^{-1}(\btheta_0)\right). 
\label{eq:MLE_normality_iid}
\end{equation}
\end{theorem}

\subsubsection{Assumptions}
\label{subsubsec:assumptions_normality_iid}
Along with the assumptions (H1) -- (H4), we further assume the following:
\begin{itemize}
\item[(H5)] The true value $\btheta_0\in int\left(\bTheta\right)$.
\item[(H6)] The Fisher's information matrix $\mathcal I(\btheta)$ is finite and non-singular, for all $\btheta\in\bTheta$.
\item[(H7)] 
Letting $\left[b'_{\bbeta}(x)\right]_k=\frac{\partial}{\partial\beta_k}b_{\bbeta}(x)$ for $k=1,\ldots,q$; 
$\left[b''_{\bbeta}(x)\right]_{kl}=\frac{\partial^2}{\partial\beta_k\partial\beta_l} b_{\bbeta}(x)$ for $k,l=1,\ldots,q$,
and $\left[b'''_{\bbeta}(x)\right]_{klm}=\frac{\partial^3}{\partial\beta_k\partial\beta_l\partial\beta_m} b_{\bbeta}(x)$ for 
$k,l,m=1,\ldots,q$, there exist constants $0<c<\infty, 0<\gamma_1, \gamma_2, \gamma_3, \gamma_4 \leq 1$ 
such that for each combination of $k,l,m=1,\ldots,q$, for any $\bbeta_1,\bbeta_2\in\mathbb R^q$, for all $x\in\mathbb R$,
$$\left|b_{\bbeta_1}(x)-b_{\bbeta_2}(x)\right|\leq c\parallel\bbeta_1-\bbeta_2\parallel^{\gamma_1};$$
$$\left|\left[b^{\prime}_{\bbeta_1}(x)\right]_k-\left[b^{\prime}_{\bbeta_2}(x)\right]_k\right|\leq c\parallel\bbeta_1-\bbeta_2\parallel^{\gamma_2};$$
$$\left|\left[b^{\prime\prime}_{\bbeta_1}(x)\right]_{kl}-\left[b^{\prime\prime}_{\bbeta_2}(x)\right]_{kl}\right|\leq c\parallel\bbeta_1-\bbeta_2\parallel^{\gamma_3};$$
$$\left|\left[b^{\prime\prime\prime}_{\bbeta_1}(x)\right]_{klm}-\left[b^{\prime\prime\prime}_{\bbeta_2}(x)\right]_{klm}\right|
\leq c\parallel\bbeta_1-\bbeta_2\parallel^{\gamma_4}.$$

\end{itemize}

\subsubsection{Verification of the above regularity conditions for asymptotic normality in our $SDE$ set-up}
\label{subsubsec:MLE_normality_iid}
In Section \ref{subsubsec:MLE_consistency_iid} almost sure consistency of the $MLE$ $\hat\btheta_n$
has been established. Hence, $\hat\btheta_n\stackrel{P}{\rightarrow}\btheta$ under $P_{\btheta}$ for all $\btheta$. 
With assumptions (H1)--(H4), (H7), Theorem B.4 of \ctn{BLSP99} and the dominated convergence theorem, 
interchangability of differentiation and integration in case of stochastic integration and usual integration 
respectively can be assured, from which it can be easily deduced that differentiation can be passed under
the integral sign, as required by Theorem \ref{theorem:theorem2}. With the same arguments, it follows that in our case 
$\frac{\partial^2}{\partial\theta_k\partial\theta_j}\log f_1(x|\btheta)$ is differentiable in 
$\btheta=(\bbeta,\xi_0)$, and the derivative has finite expectation due to compactness of the parameter space and (H7). 
Hence, (\ref{eq:h1}) and (\ref{eq:h2}) clearly hold.

In other words, asymptotic normality of the $MLE$, of the form (\ref{eq:MLE_normality_iid}), holds in our case. 
Formally,
\begin{theorem}
\label{theorem:asymp_normal_iid}
Assume the $iid$ setup and conditions (H1) -- (H7). 
Then, as $n\rightarrow\infty$, the $MLE$ is asymptotically normally distributed as
(\ref{eq:MLE_normality_iid}).
\end{theorem}
\section{Consistency and asymptotic normality of $MLE$ in the non-$iid$ set-up}
\label{sec:consistency_non_iid}

We now consider the case where the processes $X_i(\cdot);~i=1,\ldots,n$, are independently,
but not identically distributed. 
In this case, $\bxi=(\xi_{0},\xi_{1},\ldots,\xi_{p})$ where at least one of the coefficients
$\xi_{1},\ldots,\xi_{p}$ is non-zero, guaranteeing the presence of at least one time-varying covariate. 
Hence, in this set-up $\btheta=(\bbeta,\bxi)$. 

Moreover, following \ctn{Maitra14a}, \ctn{Maitra14b} we allow the initial values $x^i$ and the time limits $T_i$
to be different for $i=1,\ldots,n$, but assume that the sequences $\{T_1,T_2,\ldots\}$ and
$\{x^1,x^2,\ldots\}$ are sequences entirely contained in compact sets $\mathfrak T$ and $\mathfrak X$, respectively.
Compactness ensures that there exist convergent subsequences with limits in $\mathfrak T$ and $\mathfrak X$;
for notational convenience, we continue to denote the convergent subsequences as  $\{T_1,T_2,\ldots\}$
and $\{x^1,x^2,\ldots\}$. Thus, let the limts be $T^{\infty}\in \mathfrak T$ and $x^{\infty}\in \mathfrak X$.

Henceforth, we denote the process associated with the initial value $x$ and time point $t$ as $X(t,x)$ 
and so for $x\in \mathfrak X$ and $T\in \mathfrak T$, we let
\begin{align}
U_{\btheta}(x,T,\bz)&=\int_0^T\frac{\phi_{\bxi}b_{\bbeta}(X(s,x))}{\sigma^2(X(s,x))}d X(s,x)\label{eq:u_x_T};\\
V_{\btheta}(x,T,\bz)&=\int_0^T\frac{\phi_{\bxi}b_{\bbeta}^2(X(s,x))}{\sigma^2(X(s,x))}ds.\label{eq:v_x_T}
\end{align}
Clearly, $U_{\btheta}(x^i,T_i,\bz_i)=U_{i,{\btheta}}$ and $V_{\btheta}(x^i,T_i,\bz_i)=V_{i,{\btheta}}$, 
where $U_{i,{\btheta}}$ and $V_{i,{\btheta}}$ are given by (\ref{eq:u_v}). In this non-$iid$ set-up 
we assume, following \ctn{Maitra15}, that
\begin{itemize}
\item[(H8)]
For $l=1,\ldots,p$, and for $t\in [0,T_i]$,
\begin{equation}
\frac{1}{n}\sum_{i=1}^n g_l(z_{il}(t))\rightarrow c_{l}(t);
\label{eq:H5_prime_1}
\end{equation}
and, for $l,m=1,\ldots,p$; $t\in [0,T_i]$,
\begin{equation}
\frac{1}{n}\sum_{i=1}^n g_l(z_{il}(t))g_m(z_{im}(t))\rightarrow c_l(t)c_m(t),
\label{eq:H5_prime_2}
\end{equation}
as $n\rightarrow\infty$, where $c_l(t)$ are real constants.
\end{itemize}

For $x=x^k$, $T=T_k$ and $\bz=\bz_k$, we denote the Kullback-Leibler 
distance and the Fisher's information
as $\mathcal K_k(\btheta_0,\btheta)$ ($\mathcal K_k(\btheta,\btheta_0)$) and $\mathcal I_k(\btheta)$, respectively.
Then the following results hold in the same way as Lemma 11 of \ctn{Maitra15}.

%
%
\begin{lemma}
\label{lemma:KL_Fisher_limit}
Assume the non-$iid$ set-up, (H1) -- (H4) and (H8).
Then for any $\btheta\in\bTheta$,
\begin{align}
\underset{n\rightarrow\infty}{\lim}~\frac{\sum_{k=1}^n\mathcal K_k(\btheta_0,\btheta)}{n}&=\mathcal K(\btheta_0,\btheta);
\label{eq:kl_limit_1}\\
\underset{n\rightarrow\infty}{\lim}~\frac{\sum_{k=1}^n\mathcal K_k(\btheta,\btheta_0)}{n}&=\mathcal K(\btheta,\btheta_0);
\label{eq:kl_limit_2}\\
\underset{n\rightarrow\infty}{\lim}~\frac{\sum_{k=1}^n\mathcal I_k(\btheta)}{n}&=\mathcal I(\btheta),
\label{eq:fisher_limit_1}
\end{align}
where the limits $\mathcal K(\btheta_0,\btheta)$, $\mathcal K(\btheta,\btheta_0)$ and $\mathcal I(\btheta)$ are well-defined
Kullback-Leibler divergences and Fisher's information, respectively.
\end{lemma}

Lemma \ref{lemma:KL_Fisher_limit} will be useful in our asymptotic investigation in the non-$iid$ set-up. 
In this set-up, we first investigate consistency and asymptotic normality of $MLE$ using the results 
of \ctn{Hoadley71}. 
\subsection{Consistency and asymptotic normality of $MLE$ in the non-$iid$ set-up}
\label{subsec:consistency_non_iid}

Following \ctn{Hoadley71} we define the following:
\begin{align}
R_i(\btheta)&=\log\frac{f_i(X_i|\btheta)}{f_i(X_i|\btheta_0)}\quad\mbox{if}\ \ f_i(X_i|\btheta_0)>0\notag\\
&=0 \quad\quad\quad\quad\quad\quad\mbox{otherwise}.
\label{eq:R1}
\end{align}

\begin{align}
R_i(\btheta,\rho)&=\sup\left\{R_i(\bpsi):\|\bpsi-\btheta\|\leq\rho\right\}\label{eq:R2}\\
{\mathcal V}_i(r)&=\sup\left\{R_i(\btheta):\|\btheta\|>r\right\}.\label{eq:V}
\end{align}
Following \ctn{Hoadley71} we denote by $r_i(\btheta)$, $r_i(\btheta,\rho)$ and $v_i(r)$
to be expectations of $R_i(\btheta)$, $R_i(\btheta,\rho)$ and ${\mathcal V}_i(r)$ under $\btheta_0$; 
for any sequence $\{a_i;i=1,2,\ldots\}$ we denote
$\sum_{i=1}^na_i/n$ by $\bar a_n$.

\ctn{Hoadley71} proved that if the following regularity conditions are satisfied, then 
the MLE $\hat\btheta_n\stackrel{P}{\rightarrow}\btheta_0$:
\begin{itemize}
\item[(1)] $\bTheta$ is a closed subset of $\mathbb R^{p+q+1}$.
\item[(2)] $f_i(X_i|\btheta)$ is an upper semicontinuous 
function of $\btheta$, uniformly in $i$, 
a.s. $[P_{\btheta_0}]$.
\item[(3)] There exist $\rho^*=\rho^*(\btheta)>0$, $r>0$ and $0<K^*<\infty$ for which 
\begin{enumerate}
\item[(i)] $E_{\btheta_0}\left[R_i(\btheta,\rho)\right]^2\leq K^*,\quad 0\leq\rho\leq\rho^*$;
\item[(ii)] $E_{\btheta_0}\left[{\mathcal V}_i(r)\right]^2\leq K^*$.
\end{enumerate}
\item[(4)]
\begin{enumerate}
\item[(i)]$\underset{n\rightarrow\infty}{\lim}~\bar r_n(\btheta)<0,\quad\btheta\neq\btheta_0$;
\item[(ii)]$\underset{n\rightarrow\infty}{\lim}~\bar v_n(r)<0$.
\end{enumerate}
\item[(5)] $R_i(\btheta,\rho)$ and ${\mathcal V}_i(r)$ are measurable functions of $X_i$.
\end{itemize}
Actually, conditions (3) and (4) can be weakened but these are more easily applicable (see \ctn{Hoadley71} for details).

\subsubsection{Verification of the regularity conditions}
\label{subsubsec:consistency_non_iid}

Since $\bTheta$ is compact in our case, the first regularity condition
clearly holds. 

For the second regularity condition, note that given $X_i$, 
$f_i(X_i|\btheta)$ is continuous by our assumptions (H1) -- (H4), as already noted in 
Section \ref{subsubsec:MLE_consistency_iid}; in fact, uniformly continuous in $\btheta$ in our case, 
since $\bTheta$ is compact. Hence, for any given $\epsilon>0$, there exists $\delta_i(\epsilon)>0$, independent
of $\btheta$,
such that $\|\btheta_1-\btheta_2\|<\delta_i(\epsilon)$ implies $\left|f(X_i|\btheta_1)-f(X_i|\btheta_2)\right|<\epsilon$.
Now consider a strictly positive function $\delta_{x,T}(\epsilon)$, continuous in $x\in\mathfrak X$ and $T\in\mathfrak T$,
such that $\delta_{x^i,T_i}(\epsilon)=\delta_i(\epsilon)$. Let 
$\delta(\epsilon)=\underset{x\in\mathfrak X,T\in\mathfrak T}{\inf}\delta_{x,T}(\epsilon)$. Since
$\mathfrak X$ and $\mathfrak T$ are compact, it follows that $\delta(\epsilon)>0$. Now it holds that
$\|\btheta_1-\btheta_2\|<\delta(\epsilon)$ implies $\left|f(X_i|\btheta_1)-f(X_i|\btheta_2)\right|<\epsilon$,
for all $i$. Hence, the second regularity condition is satisfied.

Let us now focus attention on condition (3)(i).
\begin{align}
R_i(\btheta) =& U_{i,\btheta} - \frac{V_{i,\btheta}}{2} - U_{i,\btheta_0} + \frac{V_{i,\btheta_0}}{2}\notag\\
\leq & U_{i,\btheta} + \frac{V_{i,\btheta}}{2} - U_{i,\btheta_0} + \frac{V_{i,\btheta_0}}{2}.
\label{eq:upper_bound1}
\end{align}
Let us denote $\left\{\bpsi\in\mathbb R^{p+q+1}:\|\bpsi-\btheta\|\leq\rho\right\}$ by 
$S(\rho,\btheta)$. 
Here $0<\rho<\rho^*(\btheta)$, and $\rho^*(\btheta)$ 
is so small that
$S(\rho,\btheta)\subset\bTheta$ for all $\rho\in (0,\rho^*(\btheta))$. It then follows from (\ref{eq:upper_bound1}) that
\begin{align}
\underset{\psi\in S(\rho,\btheta)}{\sup}~R_i(\bpsi)
&\leq \underset{\btheta\in S(\rho,\btheta)}{\sup}~\left(U_{i,\btheta} 
+ \frac{V_{i,\btheta}}{2}\right)- U_{i,\btheta_0} + \frac{V_{i,\btheta_0}}{2}.
\label{eq:sup_R1}
\end{align}
The supremums in (\ref{eq:sup_R1}) are finite due to compactness of $S(\rho,\btheta)$. Let the supremum 
be attained at some $\btheta^*$	where $\btheta^*=\btheta^*(X_i).$ Then, the expectation of the square of the upper bound can 
be calculated in the same way as (\ref{eq:uvexp_new1}) noting that $\bar N_{\btheta}$ in this case will be $S(\rho,\btheta)$.
Since under $P_{\btheta_0}$, finiteness of moments of all orders of each term in the upper 
bound is ensured by Lemma 10 of \ctn{Maitra15}, it follows that
\begin{equation}
E_{\btheta_0}\left[R_i(\btheta,\rho)\right]^2\leq K_i(\btheta),
\label{eq:upper_bound2}
\end{equation}
where $K_i(\btheta)=K(x^i,T_i,\bz_i,\btheta)$, with $K(x,T,\bz,\btheta)$ being a continuous function of 
$(x,T,\bz,\btheta)$, continuity being again a consequence
of Lemma 10 of \ctn{Maitra15}. 
Since 
because of compactness of $\mathfrak X$, $\mathfrak T$ and $\bTheta$,
$$K_i(\btheta)\leq \underset{x\in\mathfrak X,T\in\mathfrak T,\bz\in\mathfrak Z,\btheta\in\bTheta}{\sup}
~K(x,T,\bz,\btheta)<\infty,$$
regularity condition (3)(i) follows.

To verify condition (3)(ii), first note that we can choose $r>0$ such that $\|\btheta_0\|<r$ and
$\{\btheta\in\bTheta:\|\btheta\|>r\}\neq\emptyset$. 
It then follows that 
$\underset{\left\{\btheta\in\bTheta:\|\btheta\|>r\right\}}{\sup}~R_i(\btheta)\leq \underset{\btheta\in\bTheta}{\sup}~R_i(\btheta)$
for every $i\geq 1$. The right hand side is bounded by the same expression as the right hand side of (\ref{eq:sup_R1}),
with only $S(\rho,\btheta)$ replaced with $\bTheta$. 
The rest of the verification follows in the same way as verification of (3)(i).

To verify condition (4)(i) note that by (\ref{eq:kl_limit_1}) 
\begin{equation}
\underset{n\rightarrow\infty}{\lim}~\bar r_n=-\underset{n\rightarrow\infty}{\lim}~\frac{\sum_{i=1}^n\mathcal K_i(\btheta_0,\btheta)}{n}
=-\mathcal K(\btheta_0,\btheta)<0\quad\mbox{for}~\btheta\neq\btheta_0.
\label{eq:lim1}
\end{equation}
In other words, (4)(i) is satisfied. 

Verification of (4)(ii) follows exactly in a similar way as verified in \ctn{Maitra14a} except 
that the concerned moment existence result follows from Lemma 10 of \ctn{Maitra15}. 
Regularity condition (5) is seen to hold by the same arguments as in \ctn{Maitra14a}.

In other words, in the non-$iid$ $SDE$ framework, the following
theorem holds:
\begin{theorem}
\label{theorem:consistency_non_iid}
Assume the non-$iid$ $SDE$ setup and conditions (H1) -- (H4) and (H8). 
Then it holds that $\hat\btheta_n\stackrel{P}{\rightarrow}\btheta_0$, as $n\rightarrow\infty$.
\end{theorem}

\subsection{Asymptotic normality of $MLE$ in the non-$iid$ set-up}
\label{subsec:normality_non_iid}

Let $\zeta_i(x,\btheta)=\log f_i(x|\btheta)$; also, let $\zeta'_i(x,\btheta)$ be the $(p+q+1)\times 1$ vector
with $k$-th component $\zeta'_{i,k}(x,\theta)=\frac{\partial}{\partial\theta_k}\zeta_i(x,\btheta)$, and
let $\zeta''_i(x,\btheta)$ be the $(p+q+1)\times (p+q+1)$ matrix with $(k,l)$-th element
$\zeta''_{i,kl}(x,\btheta)=\frac{\partial^2}{\partial\theta_k\partial\theta_l}\zeta_i(x,\btheta)$.

For proving asymptotic normality in the non-$iid$ framework, \ctn{Hoadley71} assumed the
following regularity conditions:
\begin{itemize}
\item[(1)] $\bTheta$ is an open subset of $\mathbb R^{p+q+1}$.
\item[(2)] $\hat\btheta_n\stackrel{P}{\rightarrow}\btheta_0$.
\item[(3)] $\zeta'_i(X_i,\btheta)$ and $\zeta''_i(X_i,\btheta)$ exist a.s. $[P_{\btheta_0}]$.
\item[(4)] $\zeta''_i(X_i,\btheta)$ is a continuous function of $\btheta$, uniformly in $i$, a.s. $[P_{\btheta_0}]$,
and is a measurable function of $X_i$.
\item[(5)] $E_{\btheta}[\zeta'_i(X_i,\btheta)]=0$ for $i=1,2,\ldots$.
\item[(6)] $\mathcal I_i(\btheta)=E_{\btheta}\left[\zeta'_i(X_i,\btheta)\zeta'_i(X_i,\btheta)^T\right]
=-E_{\btheta}\left[\zeta''_i(X_i,\btheta)\right]$, where for any vector $y$, $y^T$ denotes
the transpose of $y$.
\item[(7)] $\bar{\mathcal I}_n(\btheta)\rightarrow\bar{\mathcal I}(\btheta)$ as $n\rightarrow\infty$ and 
$\bar{\mathcal I}(\btheta)$ is positive definite.
\item[(8)] $E_{\btheta_0}\left\vert\zeta'_{i,k}(X_i,\btheta_0)\right\vert^3\leq K_2$, for some $0<K_2<\infty$.
\item[(9)] There exist $\epsilon>0$ and random variables $B_{i,kl}(X_i)$ such that
\begin{enumerate}
\item[(i)] $\sup\left\{\left\vert\zeta''_{i,kl}(X_i,\bpsi)\right\vert:\|\bpsi-\btheta_0\|\leq\epsilon\right\}
\leq B_{i,kl}(X_i)$.
\item[(ii)] $E_{\btheta_0}\left\vert B_{i,kl}(X_i)\right\vert^{1+\delta}\leq K_2$, for some $\delta>0$.
\end{enumerate}
\end{itemize}
Condition (8) can be weakened but is relatively easy to handle.
Under the above regularity conditions, \ctn{Hoadley71} prove that 
\begin{equation}
\sqrt{n}\left(\hat\btheta_n-\btheta_0\right)\stackrel{\mathcal L}{\rightarrow}
N\left(\bzero,\bar{\mathcal I}^{-1}(\btheta_0)\right).
\label{eq:MLE_normality_non_iid}
\end{equation}

\subsubsection{Validation of asymptotic normality of $MLE$ in the non-$iid$ $SDE$ set-up}
\label{subsubsec:normality_non_iid}

Condition (1) holds also for compact $\bTheta$; see \ctn{Maitra14a}. Condition (2) is a simple consequence
of Theorem \ref{theorem:consistency_non_iid}.

Conditions (3), (5) and (6) are clearly valid in our case because of interchangability 
of differentiation and integration, which follows due to (H1) -- (H4), (H7) and Theorem B.4 of \ctn{BLSP99}.
Condition (4) can be verified in exactly the same way as condition (2) of Section \ref{subsec:consistency_non_iid}
is verified; measurability of $\zeta''_i(X_i,\btheta)$ follows due to its continuity with respect to $X_i$.
Condition (7) simply follows from (\ref{eq:fisher_limit_1}). 
Compactness, continuity, and finiteness of moments
guaranteed by Lemma 10 of \ctn{Maitra15} imply conditions (8), (9)(i) and 9(ii).

In other words, in our non-$iid$ $SDE$ case we have the following theorem on asymptotic normality.
\begin{theorem}
\label{theorem:asymp_normal_non_iid}
Assume the non-$iid$ $SDE$ setup and conditions (H1) -- (H8). 
Then (\ref{eq:MLE_normality_non_iid}) holds, as $n\rightarrow\infty$. 
\end{theorem}

\section{Consistency and asymptotic normality of the Bayesian posterior in the $iid$ set-up}
\label{sec:b_consistency_iid}

\subsection{Consistency of the Bayesian posterior distribution}
\label{subsec:Bayesian_consistency_iid}

As in \ctn{Maitra14b} here we exploit Theorem 7.80 presented in \ctn{Schervish95}, stated below, to show 
posterior consistency. 
\begin{theorem}[\ctn{Schervish95}]
\label{theorem:theorem3}
Let $\{X_n\}_{n=1}^{\infty}$ be conditionally $iid$ given $\btheta$ with density $f_1(x|\btheta)$
with respect to a measure $\nu$ on a space $\left(\mathcal X^1,\mathcal B^1\right)$. Fix $\btheta_0\in\bTheta$, and define,
for each $\bM\subseteq\bTheta$ and $x\in\mathcal X^1$,
\[
Z(\bM,x)=\inf_{\bpsi\in \bM}\log\frac{f_1(x|\btheta_0)}{f_1(x|\bpsi)}.
\]
Assume that for each $\btheta\neq\btheta_0$, there is an open set $\bN_{\btheta}$ such that $\btheta\in \bN_{\btheta}$ and
that $E_{\btheta_0}Z(\bN_{\btheta},X_i)> -\infty$. 
Also assume that $f_1(x|\cdot)$ is continuous at $\btheta$ 
for every $\btheta$, a.s. $[P_{\btheta_0}]$.
For $\epsilon>0$, define 
$\bC_{\epsilon}=\{\btheta:\mathcal K_1(\btheta_0,\btheta)<\epsilon\}$, where 
\begin{equation}
\mathcal K_1(\btheta_0,\btheta)=E_{\btheta_0}\left(\log\frac{f_1(X_1|\btheta_0)}{f_1(X_1|\btheta)}\right)
\label{eq:kl1}
\end{equation}
is the Kullback-Leibler divergence measure associated with observation $X_1$.
Let $\pi$ be a prior distribution such that $\pi(\bC_{\epsilon})>0$, for every $\epsilon>0$.
Then, for every $\epsilon>0$ and open set $\bf{\mathcal N}_0$ containing $\bC_{\epsilon}$, the posterior satisfies
\begin{equation}
\lim_{n\rightarrow\infty}\pi_n\left(\bf{\mathcal N}_0|X_1,\ldots,X_n\right)=1,\quad a.s.\quad [P_{\btheta_0}].
\label{eq:posterior_consistency_iid}
\end{equation}
\end{theorem}

\subsubsection{Verification of posterior consistency}
\label{subsubsec:Bayesian_consistency_iid}
The condition $E_{\btheta_0}Z(\bN_{\btheta},X_i)> -\infty$ of the above theorem is 
verified in the context of Theorem 1 in Section \ref{subsubsec:MLE_consistency_iid}.
Continuity of the Kullback-Liebler divergence follows easily from Lemma 10 of \ctn{Maitra15}. The rest
of the verification is the same as that of \ctn{Maitra14b}.

Hence, (\ref{eq:posterior_consistency_iid}) holds in our case with any prior with positive, continuous density
with respect to the Lebesgue measure.
We summarize this result in the form of a theorem, stated below.

\begin{theorem}
\label{theorem:new_theorem3}
Assume the $iid$ set-up and conditions (H1) -- (H4). 
Let the prior distribution $\pi$ of the parameter $\theta$ satisfy $\frac{d\pi}{d\nu}=g$ 
almost everywhere on $\bTheta$, where $g(\btheta)$ is any 
positive, continuous density on $\bTheta$ with respect to the Lebesgue
measure $\nu$.  
Then the posterior (\ref{eq:posterior1}) is consistent
in the sense that for every $\epsilon>0$ and open set $\bf\mathcal N_0$ containing $\bC_{\epsilon}$, the posterior satisfies
\begin{equation}
\lim_{n\rightarrow\infty}\pi_n\left(\bf\mathcal N_0|X_1,\ldots,X_n\right)=1,\quad a.s.\quad [P_{\btheta_0}].
\label{eq:posterior_consistency_iid2}
\end{equation}
\end{theorem}

\subsection{Asymptotic normality of the Bayesian posterior distribution}
\label{subsec:Bayesian_normality_iid}
As in \ctn{Maitra14b}, we make use of Theorem 7.102 in conjunction with 
Theorem 7.89 provided in \ctn{Schervish95}. 
These theorems make use of seven regularity conditions, of which only the first four, stated below,
will be required for the $iid$ set-up. 

\subsubsection{Regularity conditions -- $iid$ case}
\label{subsubsec:regularity_iid}
\begin{itemize}
\item[(1)] The parameter space is $\bTheta\subseteq\mathbb R^{q+1}$.
\item[(2)] $\btheta_0$ is a point interior to $\bTheta$.
\item[(3)] The prior distribution of $\btheta$ has a density with respect to Lebesgue measure
that is positive and continuous at $\btheta_0$.
\item[(4)] There exists a neighborhood $\bf\mathcal N_0\subseteq\bTheta$ of $\btheta_0$ on which
$\ell_n(\btheta)= \log f(X_1,\ldots,X_n|\btheta)$ is twice continuously differentiable with 
respect to all co-ordinates of $\btheta$, 
$a.s.$ $[P_{\btheta_0}]$.
\end{itemize}

Before proceeding to justify asymptotic normality of our posterior, we furnish the relevant theorem below
(Theorem 7.102 of \ctn{Schervish95}).

\begin{theorem}[\ctn{Schervish95}]
\label{theorem:theorem4}
Let $\{X_n\}_{n=1}^{\infty}$ be conditionally $iid$ given $\btheta$.
Assume the above four regularity conditions; 
also assume that there exists $H_r(x,\btheta)$ such that, for each $\btheta_0\in int(\bTheta)$ and each
$k,j$,
\begin{align}
\sup_{\|\btheta-\btheta_0\|\leq r}\left\vert\frac{\partial^2}{\partial\theta_k\partial\theta_j}\log f_{1}(x|\btheta_0)
-\frac{\partial^2}{\partial\theta_k\partial\theta_j}\log f_{1}(x|\btheta)\right\vert\leq H_r(x,\btheta_0),
\label{eq:H1}
\end{align}
with
\begin{equation}
\lim_{r\rightarrow 0}E_{\btheta_0}H_r\left(X,\btheta_0\right)=0.
\label{eq:H2}
\end{equation}
Further suppose that
the conditions of Theorem \ref{theorem:theorem3} hold, and that the Fisher's information matrix
$\mathcal I(\btheta_0)$ is positive definite. 
Now denoting 
by $\hat\btheta_n$ the $MLE$ associated with
$n$ observations, 
let
\begin{equation}
\bSigma^{-1}_n=\left\{\begin{array}{cc}
-\ell''_n(\hat\btheta_n) & \mbox{if the inverse and}\ \ \hat\btheta_n\ \ \mbox{exist}\\
\mathbb I_{q+1} & \mbox{if not},
\end{array}\right.
\label{eq:information1}
\end{equation}
where for any $t$,
\begin{equation}
\ell''_n(t)=\left(\left(\frac{\partial^2}{\partial\theta_i\partial\theta_j}\ell_n(\btheta)\bigg\vert_{\btheta=\bt}\right)\right),
\label{eq:information2}
\end{equation}
and $\mathbb I_{q+1}$ is the identity matrix of order $q+1$.
Thus, $\bSigma^{-1}_n$ is the observed Fisher's information matrix.

Letting $\bPsi_n=\bSigma^{-1/2}_n\left(\btheta-\hat\btheta_n\right)$, it follows that for each compact subset 
$\bB$ of $\mathbb R^{q+1}$ and each $\epsilon>0$, it holds that
\begin{equation}
\lim_{n\rightarrow\infty}P_{\btheta_0}
\left(\sup_{\bPsi_n\in \bB}\left\vert\pi_n(\bPsi_n\vert X_1,\ldots,X_n)-\tilde\phi(\bPsi_n)\right\vert>\epsilon\right)=0,
\label{eq:Bayesian_normality_iid}
\end{equation}
where $\tilde\phi(\cdot)$ denotes the density of the standard normal distribution.
\end{theorem}

\subsubsection{Verification of posterior normality}
\label{subsubsec:Bayesian_normality_iid}
Observe that the 
four regularity conditions of Section \ref{subsubsec:regularity_iid} trivially hold. 
The remaining conditions of Theorem \ref{theorem:theorem4} are verified in the context of 
Theorem \ref{theorem:theorem2} in Section \ref{subsubsec:MLE_normality_iid}. 
We summarize this result in the form of the following theorem.
\begin{theorem}
\label{theorem:new_theorem4}
Assume the $iid$ set-up and conditions (H1) -- (H7). 
Let the prior distribution $\pi$ of the parameter $\theta$ satisfy $\frac{d\pi}{d\nu}=g$ 
almost everywhere on $\bTheta$, where $g(\btheta)$ is any 
density with respect to the Lebesgue measure $\nu$ which is positive and continuous at $\btheta_0$. 
Then, letting $\bPsi_n=\bSigma^{-1/2}_n\left(\btheta-\hat\btheta_n\right)$, it follows that for each compact subset 
$\bB$ of $\mathbb R^{q+1}$ and each $\epsilon>0$, it holds that
\begin{equation}
\lim_{n\rightarrow\infty}P_{\btheta_0}
\left(\sup_{\bPsi_n\in \bB}\left\vert\pi_n(\bPsi_n\vert X_1,\ldots,X_n)-\tilde\phi(\bPsi_n)\right\vert>\epsilon\right)=0.
\label{eq:Bayesian_normality_iid2}
\end{equation}
\end{theorem}

\section{Consistency and asymptotic normality of the Bayesian posterior in the non-$iid$ set-up}
\label{sec:b_consistency_non_iid}

For consistency and asymptotic normality in the non-$iid$ Bayesian framework we utilize the result 
presented in \ctn{Choi04} and Theorem 7.89 of \ctn{Schervish95}, respectively.

\subsection{Posterior consistency in the non-$iid$ set-up}
\label{subsec:posterior_consistency_non_iid}

We consider the following extra assumption for our purpose.
\begin{itemize}
\item[(H9)] There exist  strictly positive functions $\alpha_1^*(x,T,\bz,\btheta)$ and $\alpha_2^*(x,T,\bz,\btheta)$
continuous in $(x,T,\bz,\btheta)$, such that for any $(x,T,\bz,\btheta)$,
\begin{equation*}
E_{\btheta}\left[\exp\left\{\alpha_1^*(x,T,\bz,\btheta) U_{\btheta}(x,T,\bz)\right\}\right]<\infty, 
\label{eq:moment_alpha1}
\end{equation*}
and
\begin{equation*}
E_{\btheta}\left[\exp\left\{\alpha_2^*(x,T,\bz,\btheta) V_{\btheta}(x,T,\bz)\right\}\right]<\infty, 
\label{eq:moment_alpha2}
\end{equation*}
\end{itemize}

Now, let 
\begin{equation}
\alpha^*_{1,\min}=\underset{x\in\mathfrak X,T\in\mathfrak T,\bz\in\mathfrak Z,\btheta\in\bTheta}{\inf}\alpha_1^*(x,T,\bz,\btheta),
\label{eq:alpha_star1}
\end{equation}
\begin{equation}
\alpha^*_{2,\min}=\underset{x\in\mathfrak X,T\in\mathfrak T,\bz\in\mathfrak Z,\btheta\in\bTheta}{\inf}\alpha_2^*(x,T,\bz,\btheta)
\label{eq:alpha_star2}
\end{equation}
and 
\begin{equation}
\alpha=\min\left\{\alpha^*_{1,\min}, \alpha^*_{2,\min}, c^*\right\},
\label{eq:alpha2}
\end{equation}
where $0<c^*<1/16$. 

Compactness ensures that $\alpha^*_{1,\min}, \alpha^*_{2,\min}>0$, so that $0<\alpha<1/16$.
It also holds due to compactness that for $\btheta\in\bTheta$,
\begin{equation}
\underset{x\in\mathfrak X,T\in\mathfrak T,\bz\in\mathfrak Z,\btheta\in\bTheta}{\sup}
~E_{\btheta}\left[\exp\left\{\alpha U_{\btheta}(x,T,\bz)\right\}\right]<\infty. 
\label{eq:moment_alpha_1}
\end{equation}
and
\begin{equation}
\underset{x\in\mathfrak X,T\in\mathfrak T,\bz\in\mathfrak Z,\btheta\in\bTheta}{\sup}
~E_{\btheta}\left[\exp\left\{\alpha V_{\btheta}(x,T,\bz)\right\}\right]<\infty. 
\label{eq:moment_alpha_2}
\end{equation}

This choice of $\alpha$ ensuring  (\ref{eq:moment_alpha_1}) and (\ref{eq:moment_alpha_2}) will be useful in verification
of the conditions of Theorem \ref{theorem:theorem5}, which we next state.

\begin{theorem}[\ctn{Choi04}]
\label{theorem:theorem5}
Let $\{X_i\}_{i=1}^{\infty}$ be independently distributed with densities $\{f_i(\cdot|\btheta)\}_{i=1}^{\infty}$,
with respect to a common $\sigma$-finite measure, where $\btheta\in\bTheta$, a measurable space. The densities
$f_i(\cdot|\btheta)$ are assumed to be jointly measurable. Let $\btheta_0\in\bTheta$ and let $P_{\btheta_0}$
be the joint distribution of $\{X_i\}_{i=1}^{\infty}$ when $\btheta_0$ is the true value of $\btheta$.
Let $\{\bTheta_n\}_{n=1}^{\infty}$ be a sequence of subsets of $\bTheta$. Let $\btheta$ have prior $\pi$ on $\bTheta$.
Define the following:
\begin{align}
\Lambda_i(\btheta_0,\btheta) &=\log\frac{f_i(X_i|\btheta_0)}{f_i(X_i|\btheta)},\notag\\
\mathcal K_i(\btheta_0,\btheta) &= E_{\btheta_0}\left(\Lambda_i(\btheta_0,\btheta)\right)\notag\\
\varrho_i(\btheta_0,\btheta) &= Var_{\btheta_0}\left(\Lambda_i(\btheta_0,\btheta)\right).\notag
\end{align}
Make the following assumptions:
\begin{itemize}
\item[(1)] Suppose that there exists a set $\bB$ with $\pi(\bB)>0$ such that
\begin{enumerate}
\item[(i)] $\sum_{i=1}^{\infty}\frac{\varrho_i(\btheta_0,\btheta)}{i^2}<\infty,\quad\forall~\btheta\in \bB$,
\item[(ii)] For all $\epsilon>0$, $\pi\left(\bB\cap\left\{\btheta:\mathcal K_i(\btheta_0,\btheta)<\epsilon,
~\forall~i\right\}\right)>0$.
\end{enumerate}
\item[(2)] Suppose that there exist test functions $\{\Phi_n\}_{n=1}^{\infty}$, sets $\{\bOmega_n\}_{n=1}^{\infty}$
and constants $C_1,C_2,c_1,c_2>0$ such that
\begin{enumerate}
\item[(i)] $\sum_{n=1}^{\infty}E_{\btheta_0}\Phi_n<\infty$,
\item[(ii)] $\underset{\btheta\in \bTheta^c_n\cap\bOmega_n}{\sup}~E_{\btheta}\left(1-\Phi_n\right)\leq C_1e^{-c_1n}$,
\item[(iii)] $\pi\left(\bOmega^c_n\right)\leq C_2e^{-c_2n}$.
\end{enumerate}
\end{itemize}
Then,
\begin{equation}
\pi_n\left(\btheta\in \bTheta^c_n|X_1,\ldots,X_n\right)\rightarrow 0\quad a.s.~[P_{\btheta_0}].
\label{eq:posterior_consistency_non_iid}
\end{equation}
\end{theorem}

\subsubsection{Validation of posterior consistency}
\label{subsubsec:posterior_consistency_non_iid}

First note that, $f_i(X_i|\btheta)$ is given by (\ref{eq:densities}).
From the proof of Theorem \ref{theorem:consistency_non_iid}, using finiteness of moments of all orders associated with $U_{i,\theta}$ and $V_{i,\theta}$,
it follows that 
$\left|\log \frac{f_i(X_i|\btheta_0)}{f_i(X_i|\btheta)}\right|$ has an upper bound which has finite first and second order moments under $\btheta_0$, and is uniform for all 
$\btheta\in \bB$, where $\bB$ is any compact subset of $\bTheta$. Hence, for each $i$, 
$\varrho_i(\btheta_0,\btheta)$ is finite. 
Using compactness, Lemma 10 of \ctn{Maitra15} and arguments similar to that of Section 3.1.1 of \ctn{Maitra14b}, 
it easily follows that $\varrho_i(\btheta_0,\btheta)<\kappa$, for some $0<\kappa<\infty$, uniformly in $i$.
Hence, choosing a prior that gives positive probability to the set $\bB$, it follows that for all $\btheta\in \bB$,
$$ \sum_{i=1}^{\infty}\frac{\varrho_i(\btheta_0,\btheta)}{i^2}
<\kappa\sum_{i=1}^{\infty}\frac{1}{i^2}<\infty.$$ Hence, condition (1)(i) holds.
Also note that
(1)(ii) can be verified similarly as the verification of Theorem 5 of \ctn{Maitra14b}.

We now verify conditions (2)(i), (2)(ii) and (2)(iii). We let $\bOmega_n=\left(\bOmega_{1n}\times\mathbb R^{p+1}\right)$,
where $\bOmega_{1n}=\left\{\bbeta:\parallel\bbeta\parallel<M_n\right\}$, where $M_n=O(e^n)$. 
Note that
\begin{equation}
\pi\left(\bOmega^c_n\right)=\pi\left(\bOmega^c_{1n}\right)=\pi(\parallel\bbeta\parallel\geq M_n)
<E_{\pi}\left(\parallel\bbeta\parallel\right)M^{-1}_n,
\label{eq:sieve1}
\end{equation}
so that (2)(iii) holds, assuming that the prior $\pi$ is such that the expectation 
$E_{\pi}\left(\parallel\bbeta\parallel\right)$
is finite.

The verification of 2(i) can be checked in as in \ctn{Maitra14b} except the relevant changes. 
So, here we only mention the corresponding changes, skipping detailed verification.

Kolmogorov's strong law of large numbers for the non-$iid$ case (see, for example, \ctn{Serfling80}), 
holds in our problem due to finiteness of the moments of $U_{\btheta}(x,T,\bz)$ and  $V_{\btheta}(x,T,\bz)$ 
for every $x$, $T$, $\bz$ and $\btheta$ belonging to the respective compact spaces. 
Moreover, existence and boundedness of the third order derivative of 
$\ell_n(\btheta)$ with respect to its components is ensured by assumption (H7) along with compactness assumptions. 
The results stated in \ctn{Maitra15} concerned with continuity and finiteness of the moments of 
$U_{\btheta}(x,T,\bz)$ and  $V_{\btheta}(x,T,\bz)$ for every $x$, $T$, $\bz$ and $\btheta$ 
belonging to their respective compact spaces 
are needed here. The lower bound of $\log f_i(X_i|\btheta_0)-\log f_i(X_i|\hat\btheta_n)$ is denoted by 
$C_3(U_i,V_i,\hat\btheta_n)$ where 
$$C_3(U_i,V_i,\hat\btheta_n)= U_{i,\btheta_0} - \frac{V_{i,\btheta_0}}{2} - U_{i,\hat\btheta_n} - \frac{V_{i,\hat\btheta_n}}{2}$$
The rest of the verification is same as that of \ctn{Maitra14b} along with assumption (H9).

To verify condition 2(ii) we define $\bTheta_n=\bTheta_{\delta}
=\left\{(\bbeta,\bxi):\mathcal K(\btheta,\btheta_0)<\delta\right\}$,
where $\mathcal K(\btheta,\btheta_0)$, defined as in (\ref{eq:kl_limit_2}), 
is the proper Kullback-Leibler divergence. This verification is again similar to that of \ctn{Maitra14b}. 
The result can be summarized in the form of the following theorem.
\begin{theorem}
\label{theorem:new_theorem5}
Assume the non-$iid$ $SDE$ set-up. Also assume conditions (H1) -- (H9).
For any $\delta>0$, let $\bTheta_{\delta}=\left\{(\bbeta,\bxi):\mathcal K(\btheta,\btheta_0)<\delta\right\}$,
where $\mathcal K(\btheta,\btheta_0)$, defined as in (\ref{eq:kl_limit_2}), is the proper Kullback-Leibler
divergence. 
Let the prior distribution $\pi$ of the parameter $\btheta$ satisfy $\frac{d\pi}{d\nu}=h$ 
almost everywhere on $\bTheta$, where $h(\btheta)$ is any 
positive, continuous density on $\bTheta$ with respect to the Lebesgue
measure $\nu$.  
Then, as $n\rightarrow\infty$,
\begin{equation}
\pi_n\left(\btheta\in \bTheta^c_{\delta}|X_1,\ldots,X_n\right)\rightarrow 0\quad a.s.~[P_{\btheta_0}].
\label{eq:posterior_consistency_non_iid2}
\end{equation}
\end{theorem}

\subsection{Asymptotic normality of the posterior distribution in the non-$iid$ set-up} 
\label{subsec:Bayesian_normality_non_iid}

Below we present the three regularity conditions that are needed in the non-$iid$ set-up in addition
to the four conditions already stated in Section \ref{subsubsec:regularity_iid}, for
asymptotic normality given by (\ref{eq:Bayesian_normality_iid}).

\subsubsection{Extra regularity conditions in the non-$iid$ set-up}
\label{subsubsec:regularity_non_iid}

\begin{itemize}
\item[(5)] The largest eigenvalue of $\bSigma_n$ goes to zero in probability.
\item[(6)] For $\delta>0$, define $\bf\mathcal N_0(\delta)$ to be the open ball of radius $\delta$ around $\btheta_0$.
Let $\rho_n$ be the smallest eigenvalue of $\bSigma_n$. If $\bf\mathcal N_0(\delta)\subseteq\bTheta$, there exists
$K(\delta)>0$ such that
\begin{equation}
\underset{n\rightarrow\infty}{\lim}~P_{\btheta_0}\left(\underset{\btheta\in\bTheta\backslash\bf\mathcal N_0(\delta)}{\sup}~
\rho_n\left[\ell_n(\btheta)-\ell_n(\btheta_0)\right]<-K(\delta)\right)=1.
\label{eq:extra1}
\end{equation}
\item[(7)] For each $\epsilon>0$, there exists $\delta(\epsilon)>0$ such that
\begin{equation}
\underset{n\rightarrow\infty}{\lim}~P_{\btheta_0}\left(\underset{\btheta\in\bf\mathcal N_0(\delta(\epsilon)),
\|\gamma\|=1}{\sup}~
\left\vert 1+\gamma^T\bSigma^{\frac{1}{2}}_n\ell''_n(\btheta)\bSigma^{\frac{1}{2}}_n\gamma\right\vert<\epsilon\right)=1.
\label{eq:extra2}
\end{equation}
\end{itemize}
Although intuitive explanations of all the seven conditions are provided in \ctn{Schervish95}, here we briefly touch upon condition (7), which is seemingly
somewhat unwieldy. First note that condition (6) ensures consistency of the $MLE$ $\hat\btheta_n$, so that $\hat\btheta_n\in\bf\mathcal N_0(\delta)$,
as $n\rightarrow\infty$. Thus, in (7), for sufficiently large $n$ and sufficiently small $\delta(\epsilon)$, $\ell''_n(\btheta)\approx \ell''_n(\hat\btheta_n)$, 
for all $\btheta\in\bf\mathcal N_0(\delta(\epsilon))$.
Hence, from the definition of $\bSigma^{-1}_n$ given by (\ref{eq:information1}), it follows that $\bSigma^{\frac{1}{2}}_n\ell''_n(\btheta)\bSigma^{\frac{1}{2}}_n\approx-\mathbb I_{q+1}$
so that, since $\|\gamma\|=1$, $\left|1+\gamma^T\bSigma^{\frac{1}{2}}_n\ell''_n(\btheta)\bSigma^{\frac{1}{2}}_n\gamma\right|\approx \left|1-\|\gamma\|^2\right|<\epsilon$, 
for all $\btheta\in\bf\mathcal N_0(\delta(\epsilon))$ and for large enough $n$. Now it is easy to see that the role of condition (7) is only to formalize the heuristic arguments.

\subsubsection{Verification of the regularity conditions}
\label{subsubsec:verify_last}

Assumptions (H1) --  (H9), along with Kolmogorov's strong law of large numbers, are sufficient 
for the regularity conditions to hold; the arguments remain similar as those in Section 3.2.2 of \ctn{Maitra14b}. 
We provide our result in the form of the following theorem.
\begin{theorem}
\label{theorem:new_theorem6}
Assume the non-$iid$ set-up and conditions (H1) -- (H9). Let the prior distribution $\pi$ of the parameter 
$\btheta$ satisfy $\frac{d\pi}{d\nu}=h$ 
almost everywhere on $\bTheta$, where $h(\btheta)$ is any 
density with respect to the Lebesgue measure $\nu$ which is positive and continuous at $\btheta_0$. 
Then, letting $\bPsi_n=\bSigma^{-1/2}_n\left(\btheta-\hat\btheta_n\right)$, for each compact subset 
$\bB$ of $\mathbb R^{p+q+1}$ and each $\epsilon>0$, the following holds:
\begin{equation}
\lim_{n\rightarrow\infty}P_{\btheta_0}
\left(\sup_{\bPsi_n\in \bB}\left\vert\pi_n(\bPsi_n\vert X_1,\ldots,X_n)-\tilde\phi(\bPsi_n)\right\vert>\epsilon\right)=0.
\label{eq:Bayesian_normality_non_iid2}
\end{equation}
\end{theorem}

\section{Random effects $SDE$ model}
\label{sec:another}
We now consider the following system of $SDE$ models for $i=1,2,\ldots,n$:
\begin{equation}
d X_i(t)=\phi_{\bxi^i}(t)b_{\bbeta}(X_i(t))dt+\sigma(X_i(t))dW_i(t)
\label{eq:sdenew}
\end{equation}
Note that this model is the same as described in Section \ref{sec:model} except that the
parameters $\bxi^i$ now depend upon $i$. 
Indeed, now $\phi_{\bxi^i}(t)$ is given by 
\begin{equation}
\phi_{\bxi^i}(t)=\phi_{\bxi^i}(\bz_i(t))
=\xi^i_0+\xi_{1}^ig_1(z_{i1}(t))+\xi_{2}^ig_2(z_{i2}(t))+\cdots+\xi_{p}^ig_p(z_{ip}(t)),
\label{eq:phi_new}
\end{equation}
where $\bxi^i=(\xi^i_0,\xi_{1}^i,\ldots,\xi_{p}^i)^T$ is the random effect corresponding to the 
$i$-th individual for $i=1,\ldots,n$, 
and $\bz_i(t)$ is the same as in Section \ref{subsec:covariates}. 
We let $b_{\bbeta}(X_i(t),\phi_{\bxi^i})=\phi_{\bxi^i}(t)b_{\bbeta}(X_i(t))$. 
Note that our likelihood is the product over $i=1,\ldots,n$, of the following individual densities:
$$f_{i,\bxi^i,\bbeta}(X_i)=\exp\left(U_{i,\bxi^i,\bbeta}-\frac{V_{i,\bxi^i,\bbeta}}{2}\right),$$
where
$$U_{i,\bxi^i,\bbeta}=\int_0^{T_i}\frac{\phi_{\bxi^i}(s)b_{\bbeta}(X_i(s))}{\sigma^2(X_i(s))}dX_i(s)
\quad\quad\mbox{and}\quad\quad V_{i,\bxi^i,\bbeta}=
\int_0^{T_i}\frac{\phi_{\bxi^i}^2(s)b^2_{\bbeta}(X_i(s))}{\sigma^2(X_i(s))}ds.$$

Now, let $m^{\bbeta}(\bz(t),x(t))=(m^{\bbeta}_0,m^{\bbeta}_1(z_1(t),x(t)),\ldots,m^{\bbeta}_p(z_p(t),x(t)))^T$ 
be a function from 
$\mathfrak Z\times\mathbb R\rightarrow \mathbb R^{p+1}$ 
where $m^{\bbeta}_0\equiv 1$ and $m^{\bbeta}_k(z(t),x(t))=g_k(z_{k}(t))b_{\bbeta}(x(t))$; $k=1,\ldots,p$.
With this notation, the likelihood can be re-written as the product over $i=1,\ldots,n$, of the following:
\begin{equation}
f_{i,\bxi^i,\bbeta}(X_i)=\exp((\bxi^i)^T\bA^{\bbeta}_i-(\bxi^i)^T\bB^{\bbeta}_i\bxi^i) 
\label{eq:liknew}
\end{equation}
where 
\begin{equation}
\bA^{\bbeta}_i=\int_0^{T_i}\frac{m^{\bbeta}(z(s),X_i(s))}{\sigma^2(X_i(s))}dX_i(s)
\label{eq:suff3}
\end{equation}
and
\begin{equation}
\bB^{\bbeta}_i=\int_0^{T_i}\frac{m^{\bbeta}(z(s),X_i(s))\left(m^{\bbeta}\right)^T(z(s),X_i(s))}{\sigma^2(X_i(s))}ds
\label{eq:suff4}
\end{equation}
are $(p+1)\times 1$ random vectors and positive definite $(p+1)\times (p+1)$ random matrices respectively.

We assume that $\bxi^i$ are $iid$ Gaussian vectors, with expectation vector $\bmu$
and covariance matrix $\bSigma\in \bS_{p+1} (\mathbb R)$ where $\bS_{p+1} (\mathbb R)$ 
is the set of real positive definite symmetric matrices of order $p+1$. The parameter set
is denoted by $\btheta=(\bmu,\bSigma,\bbeta)\in \bTheta\subset\mathbb R^{p+1} \times\bS_{p+1}(\mathbb R)\times\mathbb R^q$.

To obtain the likelihood involving $\btheta$ we refer to the multidimensional random effects set-up of \ctn{Maud12}. 
Following Lemma 2 of \ctn{Maud12} it then follows in our case that, 
for each $i\geq 1$ and for all $\btheta$, $\bB^{\bbeta}_i+\bSigma^{-1},
\mathbb I_{p+1}+\bB^{\bbeta}_i\bSigma,\mathbb I_{p+1}+\bSigma \bB^{\bbeta}_i$ 
are invertible.

Setting $\left(\bR^{\bbeta}_i\right)^{-1}=(\mathbb I_{p+1}+\bB^{\bbeta}_i\bSigma)^{-1}\bB^{\bbeta}_i$ we obtain
\begin{align}
f_i(X_i|\btheta)&=\frac{1}{\sqrt{\det(\mathbb I_{p+1}+\bB^{\bbeta}_i\bSigma)}}
\exp\left(-\frac{1}{2}\left(\bmu-\left(\bB^{\bbeta}_i\right)^{-1}\bA^{\bbeta}_i\right)^T
\left(\bR^{\bbeta}_i\right)^{-1}\left(\bmu-\left(\bB^{\bbeta}_i\right)^{-1}\bA^{\bbeta}_i\right)\right)\notag\\
&\qquad\qquad\times\exp\left(\frac{1}{2}\left(\bA^{\bbeta}_i\right)^T\left(\bB^{\bbeta}_i\right)^{-1}\bA^{\bbeta}_i\right)
\label{eq:likelihoodnew1}
\end{align}
as our desired likelihood after integrating (\ref{eq:liknew}) with respect to the distrbution of $\bxi ^i$.

With reference to \ctn{Maud12} in our case 
$$\gamma_i(\btheta)=(\mathbb I_{p+1}+\bSigma \bB^{\bbeta}_i)^{-1}(\bA^{\bbeta}_i-\bB^{\bbeta}_i\bmu)
\quad\quad\mbox{and}\quad I_i(\bSigma)=(\mathbb I_{p+1}+\bSigma \bB^{\bbeta}_i)^{-1}\bB^{\bbeta}_i.$$
Hence, Proposition (10)(i) of \ctn{Maud12} can be seen to be hold here in a similar way by replacing 
$U_i$ and $V_i$ by $\bA^{\bbeta}_i$ and $\bB^{\bbeta}_i$ respectively.

Asymptotic investigation regarding consistency and asymptotic normality of $MLE$ and 
Bayesian posterior consistency and asymptotic posterior normality in both $iid$ and non-$iid$ set-ups 
can be established as in the one dimensional cases in \ctn{Maitra14a} and \ctn{Maitra14b} 
with proper multivariate modifications by replacing $U_i$ and $V_i$ with $\bA^{\bbeta}_i$ and $\bB^{\bbeta}_i$ respectively,
and exploiting assumptions (H1) -- (H9).

\section{Simulation studies}
\label{sec:simulated_data}

We now supplement our asymptotic theory with simulation studies where the data is generated from a specific system of $SDE$s with one covariate, with given values of the parameters. 
Specifically, in the classical case, we obtain the distribution of the $MLE$s using parametric bootstrap, along with the 95\% confidence intervals of the parameters. We demonstrate in particular
that the true values of the parameters are well-captured by the respective 95\% confidence intervals. In the Bayesian counterpart, we obtain the posterior distributions of the parameters
along with the respective 95\% credible intervals, and show that the true values fall well within the respective 95\% credible intervals.

\subsection{Distribution of $MLE$ when $n=20$}
\label{subsec:C_N20}

To demonstrate the finite sample analogue of asymptotic distribution of $MLE$  as $n\rightarrow\infty$, we consider $n=20$ individuals, where the $i$-th one is modeled by
\begin{equation}
dX_i(t)=(\theta_1+\theta_2z_{i1}(t))(\theta_3+\theta_4X_i(t))dt+\sigma dW_i(t),
\label{eq:sde2_appl}
\end{equation}
for $i=1,\ldots,20$. We fix our diffusion coefficient as $\sigma=1$.  We consider the initial value $X(0)=0$ and the time interval $[0, T]$ with $T=1$. 
Further, we choose the true values as $\theta_1=1, \theta_2= -1, \theta_3=2, \theta_4= -2$.

We assume that the time dependent covariates $z_{i1}(t)$ satisfy the following $SDE$
\begin{equation}
dz_{i1}(t)=\xi_{i1}z_{i1}(t))dt+ dW_i(t),
\label{eq:covariate_appl}
\end{equation}
for $i=1,\ldots,20$, where the coeffiicients $\xi_{i1}\stackrel{iid}{\sim} N(7,1)$ for $i=1,\ldots,20$. After simulating the covariates using the system of $SDE$s (\ref{eq:covariate_appl}), 
we generate the data using the system of $SDE$s (\ref{eq:sde2_appl}). In both the cases we discretize the time interval $[0,1]$ into $100$ equispaced time points.

The distributions of the $MLE$s of the four parameters are obtained through the parametric bootstrap method. In this method we simulated the data 1000 times by simulating
as many paths of the Brownian motion, where each data set consists of 20 individuals. Under each data set we perform the ``block-relaxation" method (see, for example,
\ctn{Lange10} and the references therein) to obtain the $MLE$. In a nutshell,
starting with some sensible initial value belonging to the parameter space, the block-relaxation method iteratively maximizes the optimizing function (here, the log-likelihood), 
successively, with respect to one parameter, fixing the others at their current values, until convergence is attained with respect to the iterations. Details follow.

For the initial values of the $\theta_j$ for $j=1,\ldots, 4$, required to begin the block-relaxation method, we simulate four $N(0,1)$ variates independently, and set them 
as the initial values $\theta^{(0)}_j$; $j=1,\ldots,4$. Denoting by $\theta^{(i)}_j$ the value of $\theta_j$ at the $i$-th iteration, for $i\geq 1$, and letting $L$ be the
likelihood, the block-relaxation method consists of the following steps:
\begin{algorithm}
\caption{Block-relaxation for $MLE$ in $SDE$ system with covariates}
\label{algo:block}
\begin{itemize}
\item[(1)] At the $i$-th iteration, for $j=1,2,3,4$, obtain $\theta^{(i)}_j$ by solving the equation $\frac{\partial \log L}{\partial \theta_j}=0$, conditionally on 
$\theta_1=\theta^{(i)}_1, \theta_2=\theta^{(i)}_2,\ldots,\theta_{j-1}=\theta^{(i)}_{j-1},\theta_{j+1}=\theta^{(i-1)}_{j+1},\ldots,\theta_4=\theta^{(i-1)}_4.$
Let $\btheta^{(i)}=\left(\theta^{(i)}_1,\ldots,\theta^{(i)}_4\right)$.
\item[(2)] Letting $\|\cdot\|$ denote the Euclidean norm, if $\left\|\btheta^{(i)}-\btheta^{(i-1)}\right\|\leq 10^{-5}$, set $\hat\btheta=\btheta^{(i)}$, where
$\hat\btheta$ stands for the maximum likelihood estimator of $\btheta=\left(\theta_1,\ldots,\theta_4\right)$.
\item[(3)] If, on the other hand, $\left\|\btheta^{(i)}-\btheta^{(i-1)}\right\|> 10^{-5}$, increase $i$ to $i+1$ and continue steps (1) and (2). 
\end{itemize}
\end{algorithm}

Once we obtain the $MLE$ by the above block-relaxation algorithm, we then plug-in $\btheta=\hat\btheta$ in (\ref{eq:sde2_appl}) and generate $1000$ data sets
from the resulting system of $SDE$s, and apply the block-relaxation algorithm to each such data set to obtain the $MLE$ associated with the data sets.
Thus, we obtain the distribution of the $MLE$ using the parametric bootstrap method. 

\begin{figure}
\begin{subfigure}
\centering
\includegraphics[height=8cm,width=7cm]{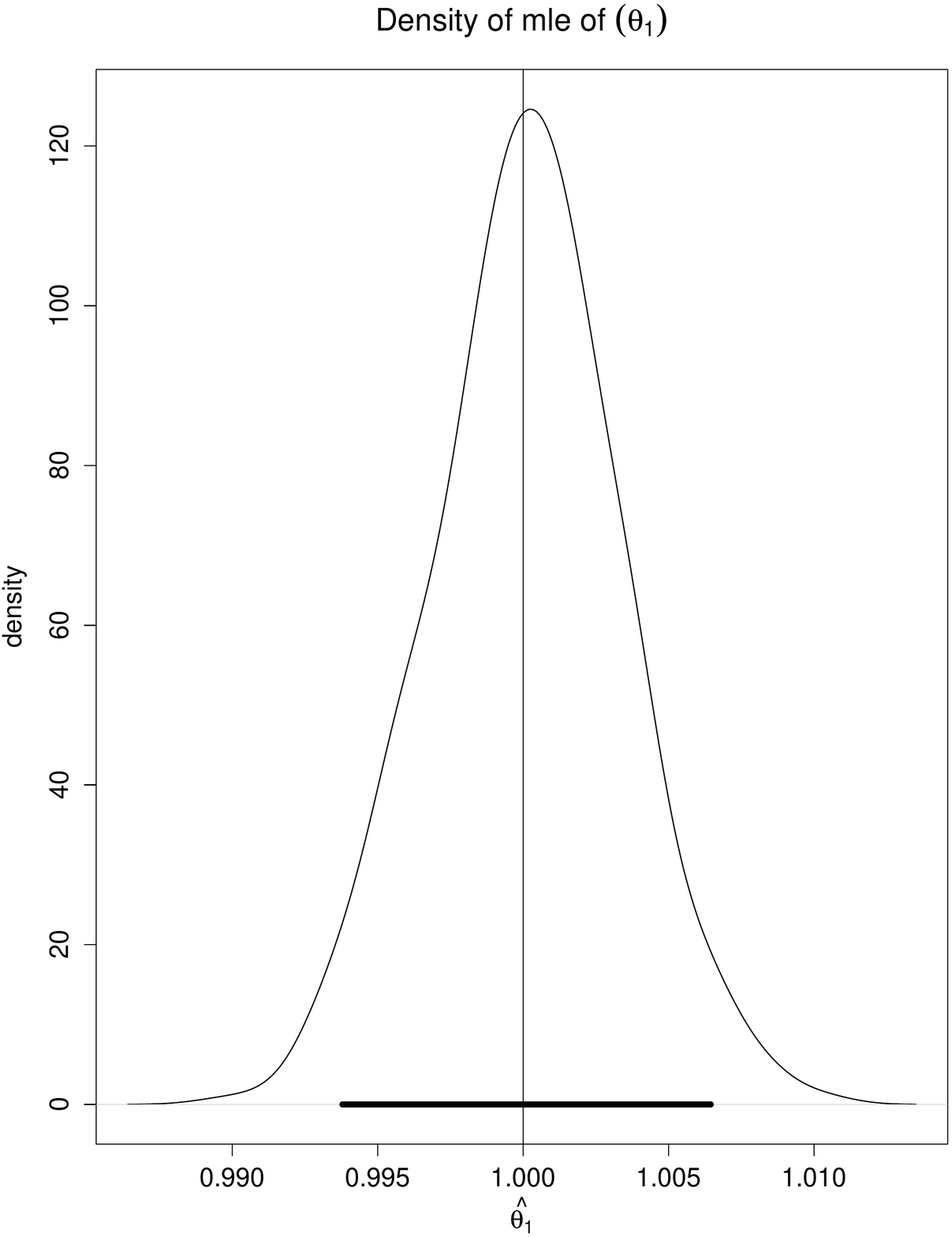}
\label{fig:sfig1}
\end{subfigure}
\begin{subfigure}
\centering
\includegraphics[height=8cm,width=7cm]{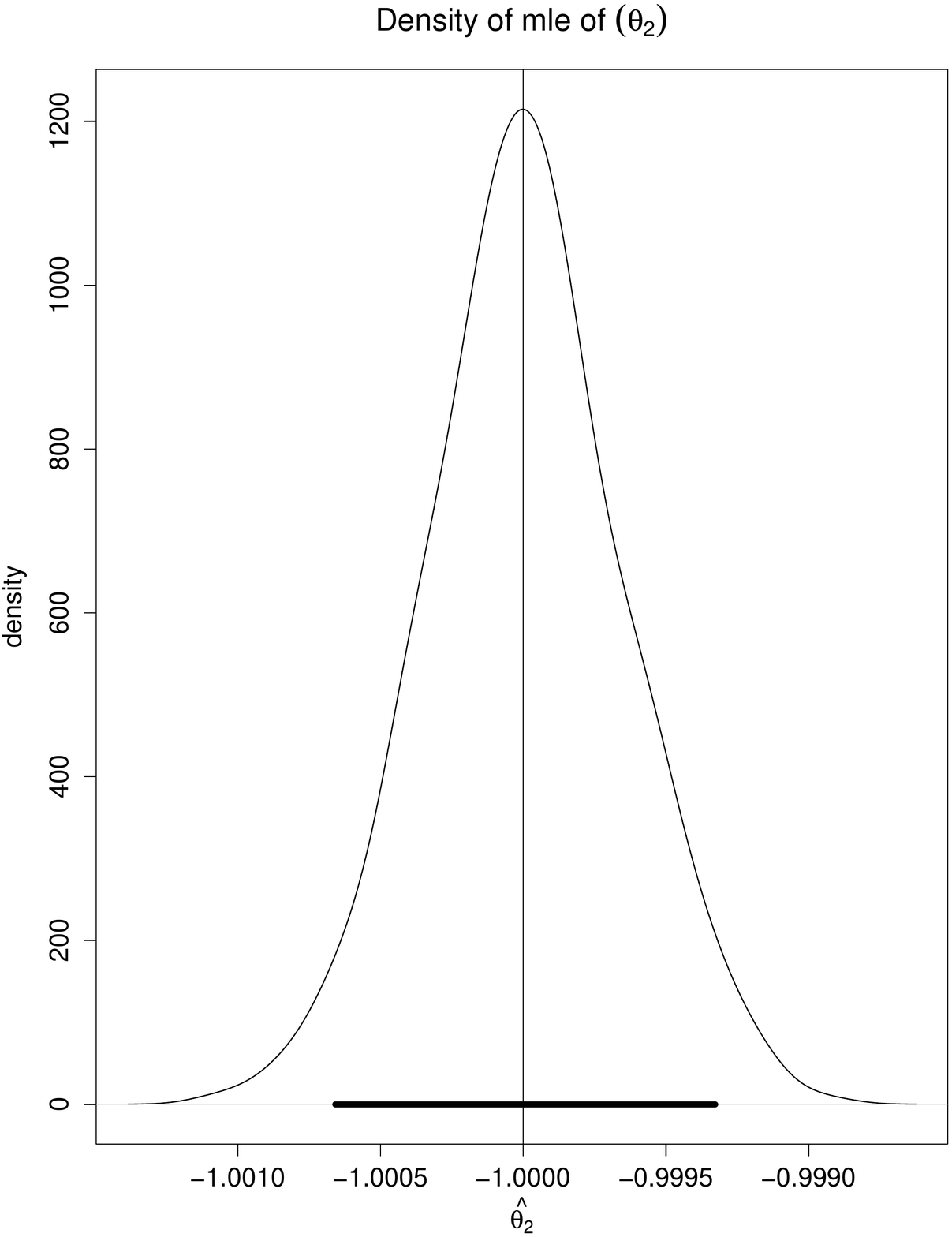}
\label{fig:sfig2}
\end{subfigure}
\begin{subfigure}
\centering
\includegraphics[height=8cm,width=7cm]{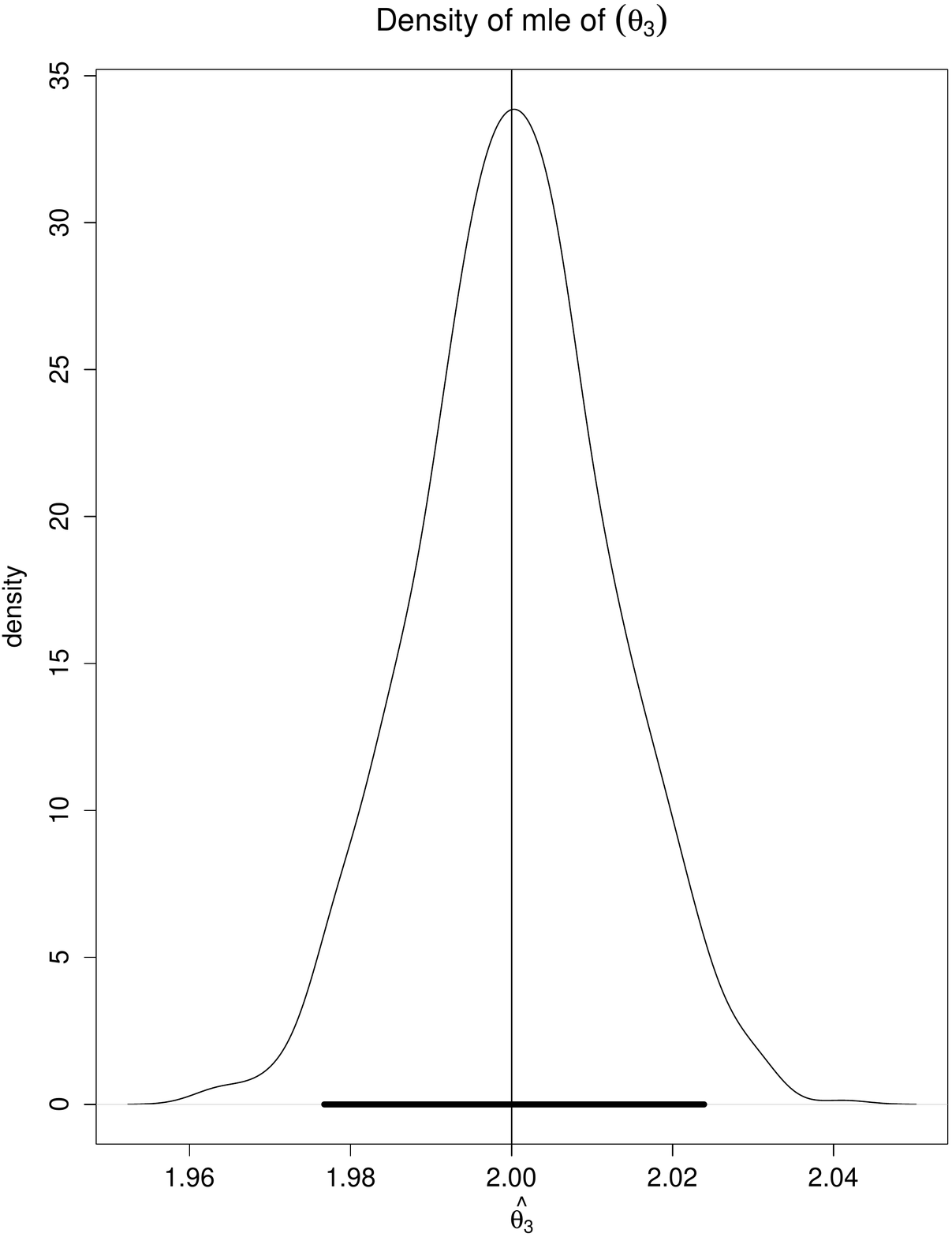}
\label{fig:sfig3}
\end{subfigure}
\begin{subfigure}
\centering
\includegraphics[height=8cm,width=7cm]{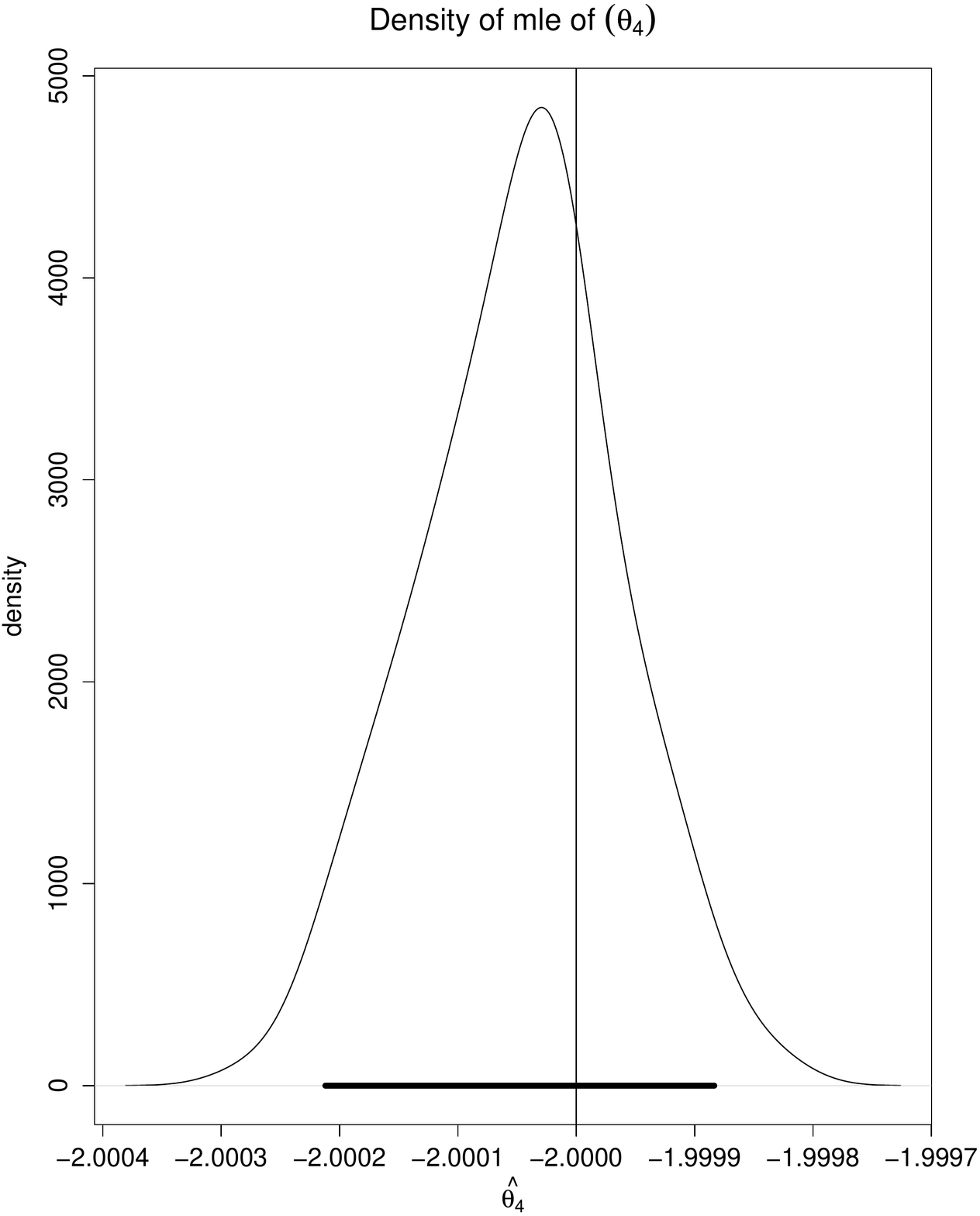}
\label{fig:sfig4}
\end{subfigure}
\label{fig:sim_c_20}
\caption{Distributions of the $MLE$s.}
\end{figure}

The distributions of the components of $\hat\btheta$ (denoted by $\hat\theta_i$ for $i=1,\ldots, 4$) are shown in Figure \ref{fig:sim_c_20}, where the associated $95\%$ confidence intervals 
are shown in bold lines. As exhibited by the figures, the 95\% confidence intervals clearly contain the true values of the respective components of $\btheta$.

\subsection{Posterior distribution of the parameters when $n=20$}
\label{subsec:B_N20}

We now consider simulation study for the Bayesian counterpart, using the same data set simulated from the system of $SDE$s given by (\ref{eq:sde2_appl}), with same
covariates simulated from (\ref{eq:covariate_appl}). We consider an empirical Bayes prior based on the $MLE$ such that for $j=1,\ldots,4$, 
$\theta_j\sim N\left(\hat\theta_j,\hat\sigma^2_j\right)$ independently,
where $\hat\theta_j$ is the $MLE$ of $\theta_j$ and $\hat\sigma^2_j$ is such that the length of the 95\% confidence interval associated with the distribution of the $MLE$ $\hat\theta_j$,
after adding one unit to both lower and upper ends of the interval, is the same as the length of the 95\% prior credible interval $[\hat\theta_j-1.96\hat\sigma_j-1,\hat\theta_j+1.96\hat\sigma_j+1]$. 
In other words, we select $\hat\sigma_j$ such that the length of the corresponding 95\% prior credible interval is the same as that of the enlarged 95\% confidence interval associated
with the distribution  of the corresponding $MLE$.
%

\begin{figure}
\begin{subfigure}
\centering
\includegraphics[height=8cm,width=7cm]{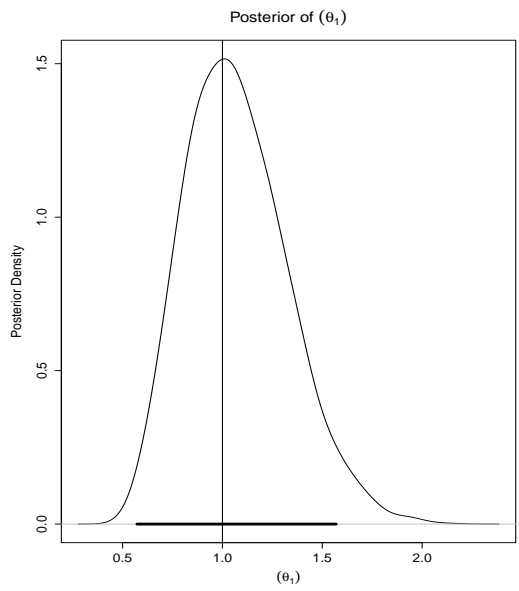}
\label{fig:sfig5}
\end{subfigure}
\begin{subfigure}
\centering
\includegraphics[height=8cm,width=7cm]{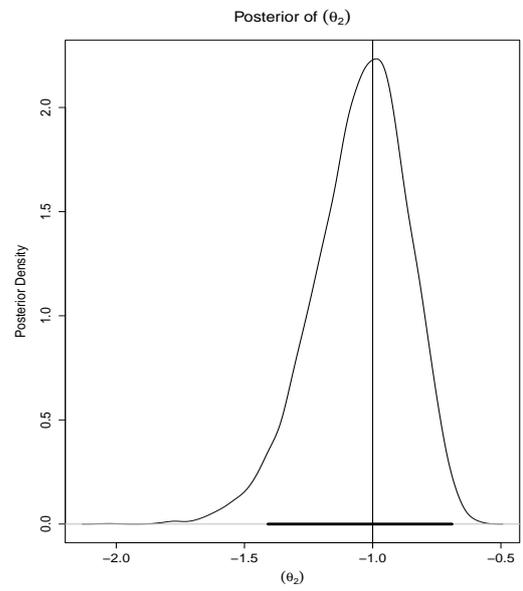}
\label{fig:sfig6}
\end{subfigure}
\begin{subfigure}
\centering
\includegraphics[height=8cm,width=7cm]{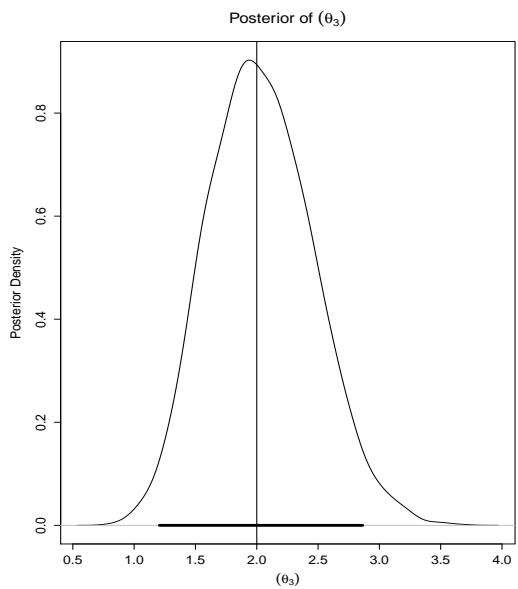}
\label{fig:sfig7}
\end{subfigure}
\begin{subfigure}
\centering
\includegraphics[height=8cm,width=7cm]{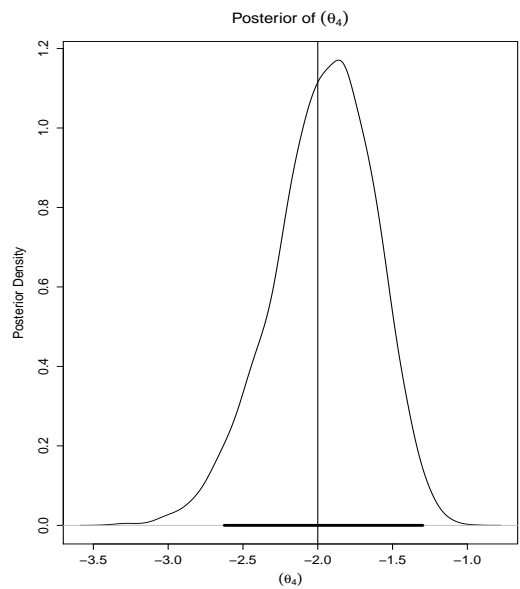}
\label{fig:sfig8}
\end{subfigure}
\label{fig:sim_b_20}
\caption{Posterior distributions of the parameters.}
\end{figure}

To simulate from the posterior distribution of $\btheta$, we perform approximate Bayesian computation (ABC) (\ctn{Tavare97}, \ctn{Beaumont02}, \ctn{marjoram03}), since the standard Markov chain
Monte Carlo (MCMC) based simulation techniques, such as Gibbs sampling and Metropolis-Hastings algorithms (see, for example, \ctn{Robert04}, \ctn{Brooks11}) 
failed to ensure good mixing behaviour of the underlying Markov chains.  
Denoting the true data set by $X_{true}$, our method of ABC is described by following steps.
\begin{algorithm}
	\caption{ABC for $SDE$ system with covariates}
	\label{algo:ABC}
\begin{itemize}
\item[(1)] For $j=1,\ldots,4$, we simulate the  parameters $\theta_j$ from their respective prior distributions.
\item[(2)] With the obtained values of the parameters we simulate the new data, which we denote by $X_{new}$, using the system of $SDE$s (\ref{eq:sde2_appl}).  
\item[(3)] We calculate the average Euclidean distance between $X_{new}$ and $X_{true}$ and denote it by $d_x$.
\item[(4)] Until $d_x<0.1$, we repeat steps (1)--(3). 
\item[(5)] Once $d_x<0.1$, we set the corresponding $\btheta$ as a realization from the posterior of $\btheta$ with approximation error $0.1$.
\item[(6)] We obtain $10000$ posterior realizations of $\btheta$ by repeating steps (1)--(5).
\end{itemize}
\end{algorithm}
Figure \ref{fig:sim_b_20} shows the posterior distribution of the parameters $\theta_i$, for $i=1,\ldots, 4$, where the $95\%$ posterior credible intervals are shown in bold lines.
Observe that all the true values of $\theta_j$; $j=1,\ldots, 4$, fall comfortably within the respective 95\% posterior credible intervals.

\section{Application to real data}
\label{sec:truedata} 

We now consider application of our $SDE$ system consisting of covariates to a real, stock market data ($467$ observations from August $5$, 2013, to June $30$, 2015) for $15$ companies. 
The data are available at {\it www.nseindia.com}.

Each company-wise data is modeled by the availabe standard financial $SDE$ models with the ``fitsde" package in $R$. 
The minimum value of BIC (Bayesian Information Criterion) is found corresponding to the $CKLS$ (Chan, Karolyi, Longstaff and Sander; see \ctn{CHan92}) model. 
Denoting the data by $X(t)$, the $CKLS$ model is described by
$$dX(t)=(\theta_1+\theta_2X(t))dt+\theta_3X(t)^{\theta_4}dW(t).$$
In our application we treat the diffusion coefficient as a fixed quantity. 
So, we fix the values of $\theta_3$ and $\theta_4$ as obtained by the ``fitsde" function, We denote $\theta_3=A$, $\theta_4=B$.

We consider the ``close price" of each company as our data $X(t)$.  IIP general index, bank interest rate and US dollar exchange rate are considered as time dependent covariates which we 
incorporate in the $CKLS$ model.

The three covariates are denoted by $c_1,c_2,c_3$, respectively. Now, our considered system of $SDE$ models
for national stock exchange data associated with the $15$ companies is the
following:
\begin{equation}
dX_i(t)=(\theta_1+\theta_2c_1(t)+\theta_3c_2(t)+\theta_4c_3(t))(\theta_5+\theta_6X_i(t))dt+A^iX_i(t)^{B^i}dW_i(t),
\label{eq:sde4_appl}
\end{equation}
for $i=1,\ldots,15$. 

\subsection{Distribution of $MLE$}
\label{subsec:C_R15}

\begin{figure}
\begin{subfigure}
\centering
\includegraphics[height=6cm,width=5cm]{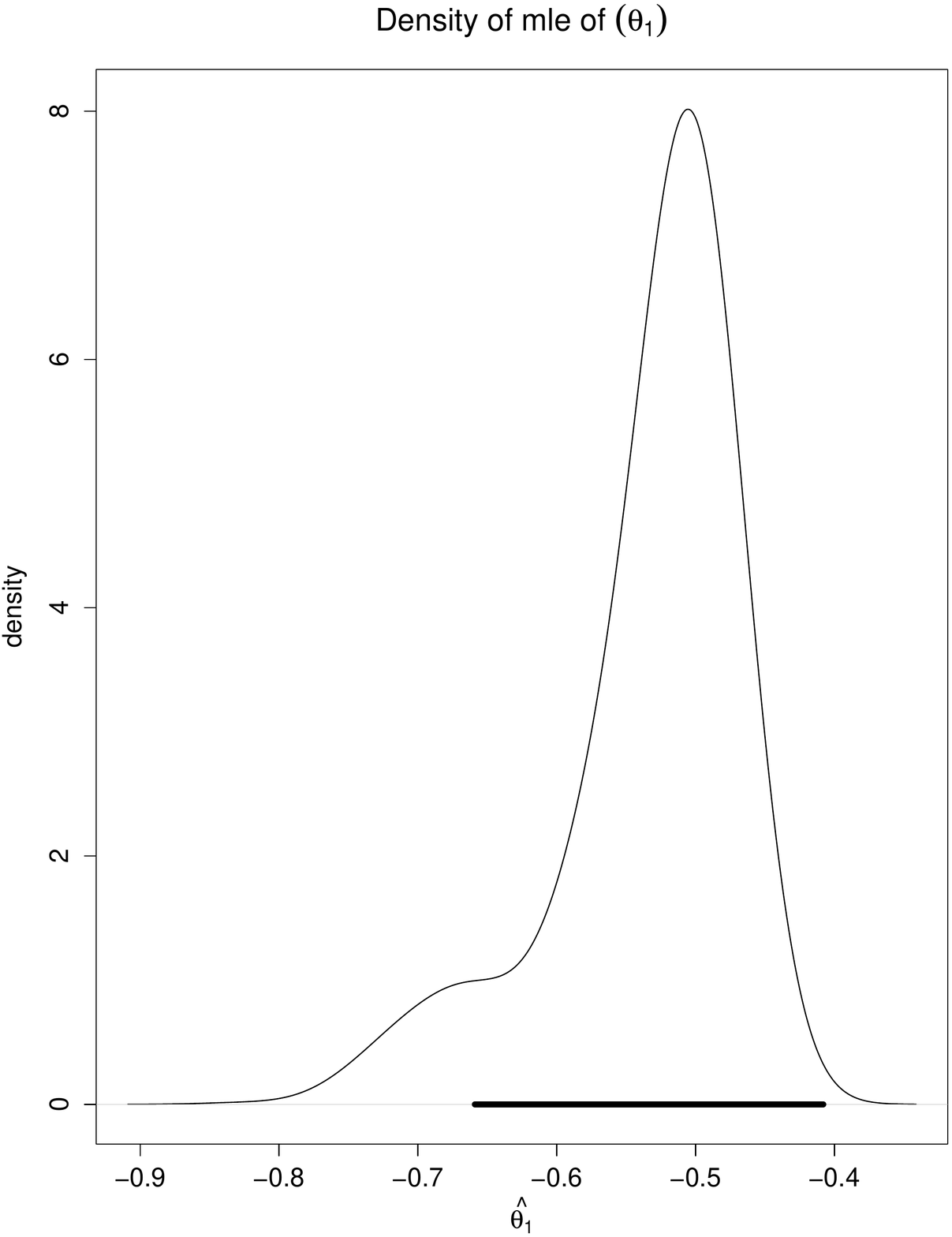}
\label{fig:sfig_r1}
\end{subfigure}
\begin{subfigure}
\centering
\includegraphics[height=6cm,width=5cm]{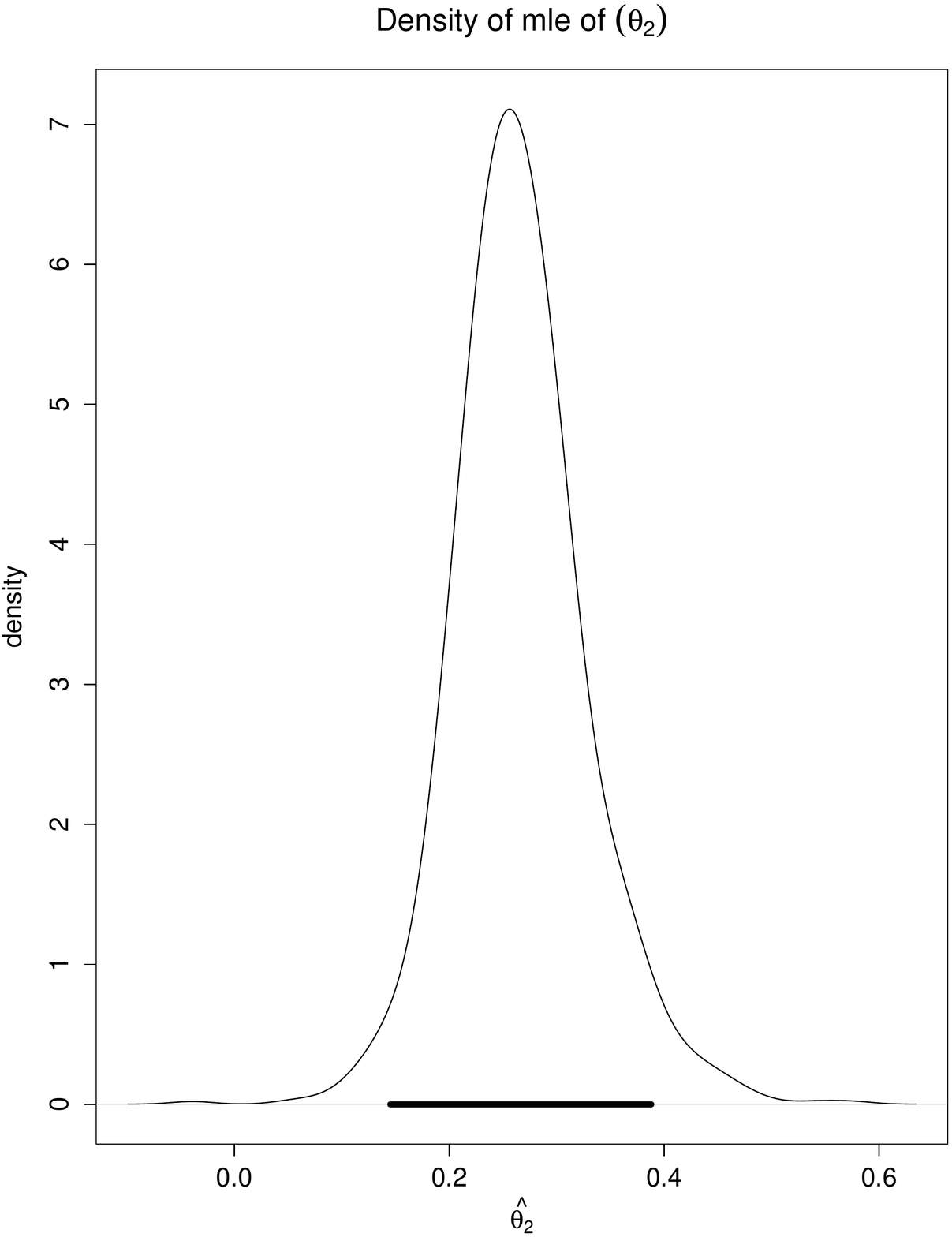}
\label{fig:sfig_r2}
\end{subfigure}
\begin{subfigure}
\centering
\includegraphics[height=6cm,width=5cm]{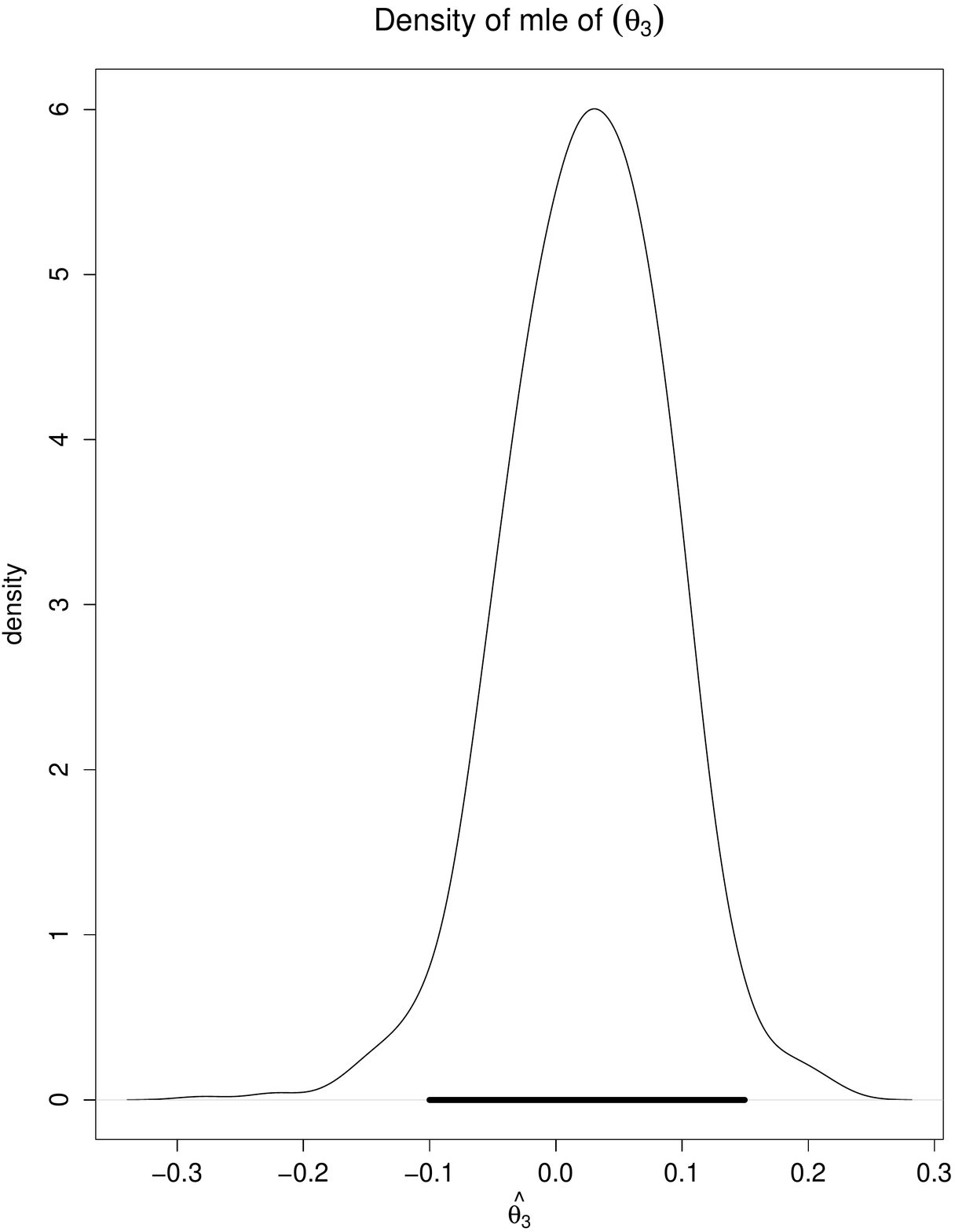}
\label{fig:sfig_r3}
\end{subfigure}
\begin{subfigure}
\centering
\includegraphics[height=6cm,width=5cm]{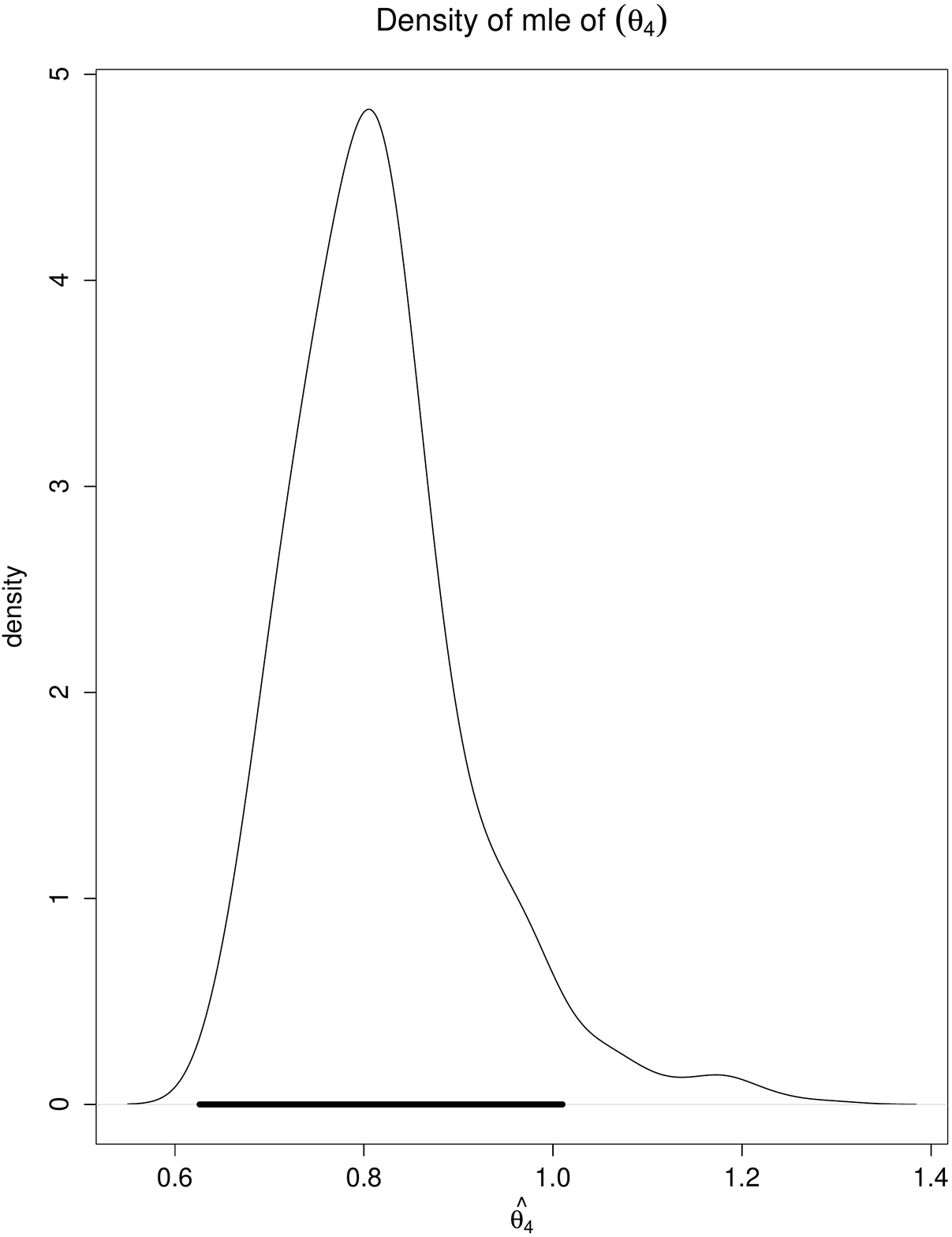}
\label{fig:sfig_r4}
\end{subfigure}
\begin{subfigure}
\centering
\includegraphics[height=6cm,width=5cm]{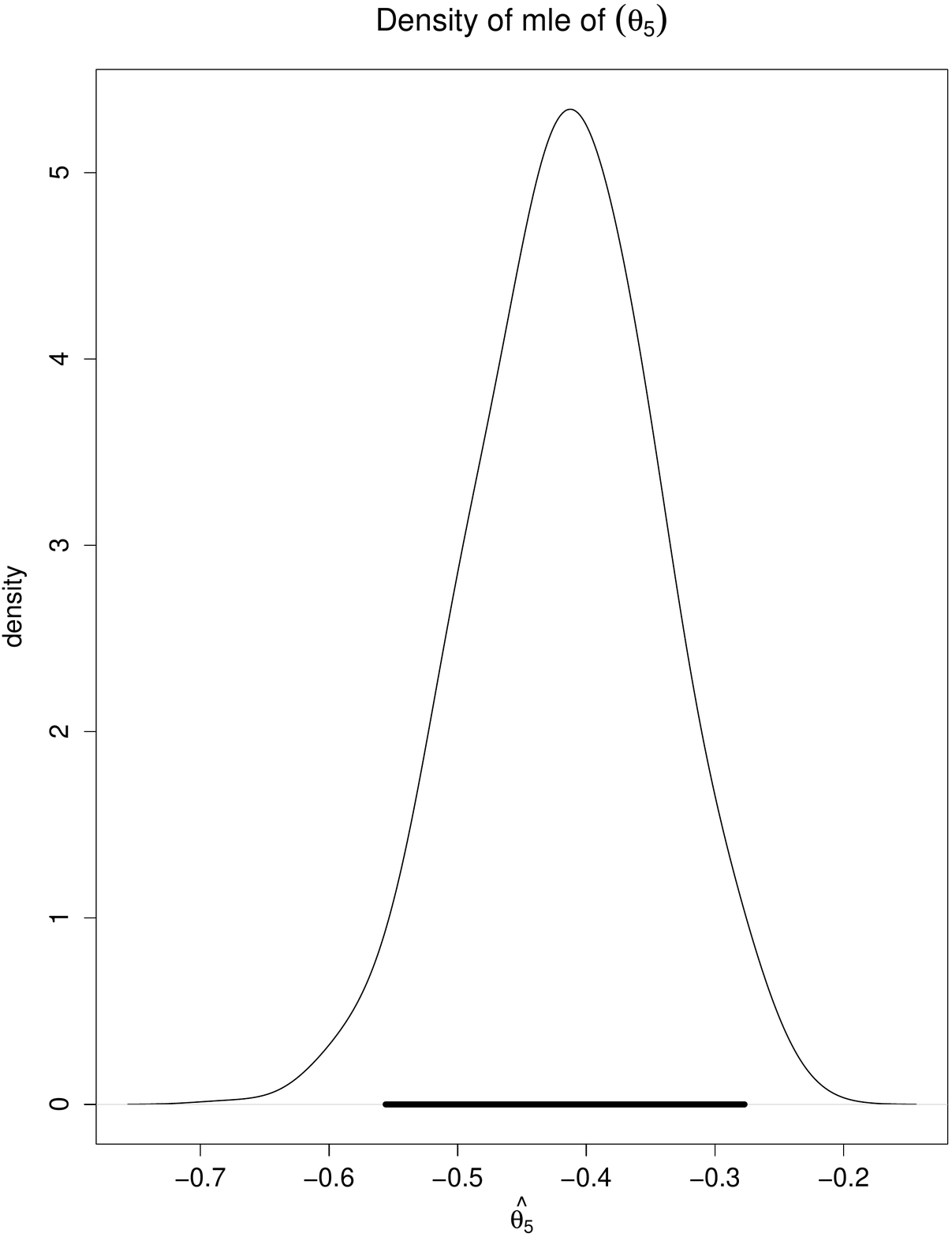}
\label{fig:sfig_r5}
\end{subfigure}
\begin{subfigure}
\centering
\includegraphics[height=6cm,width=5cm]{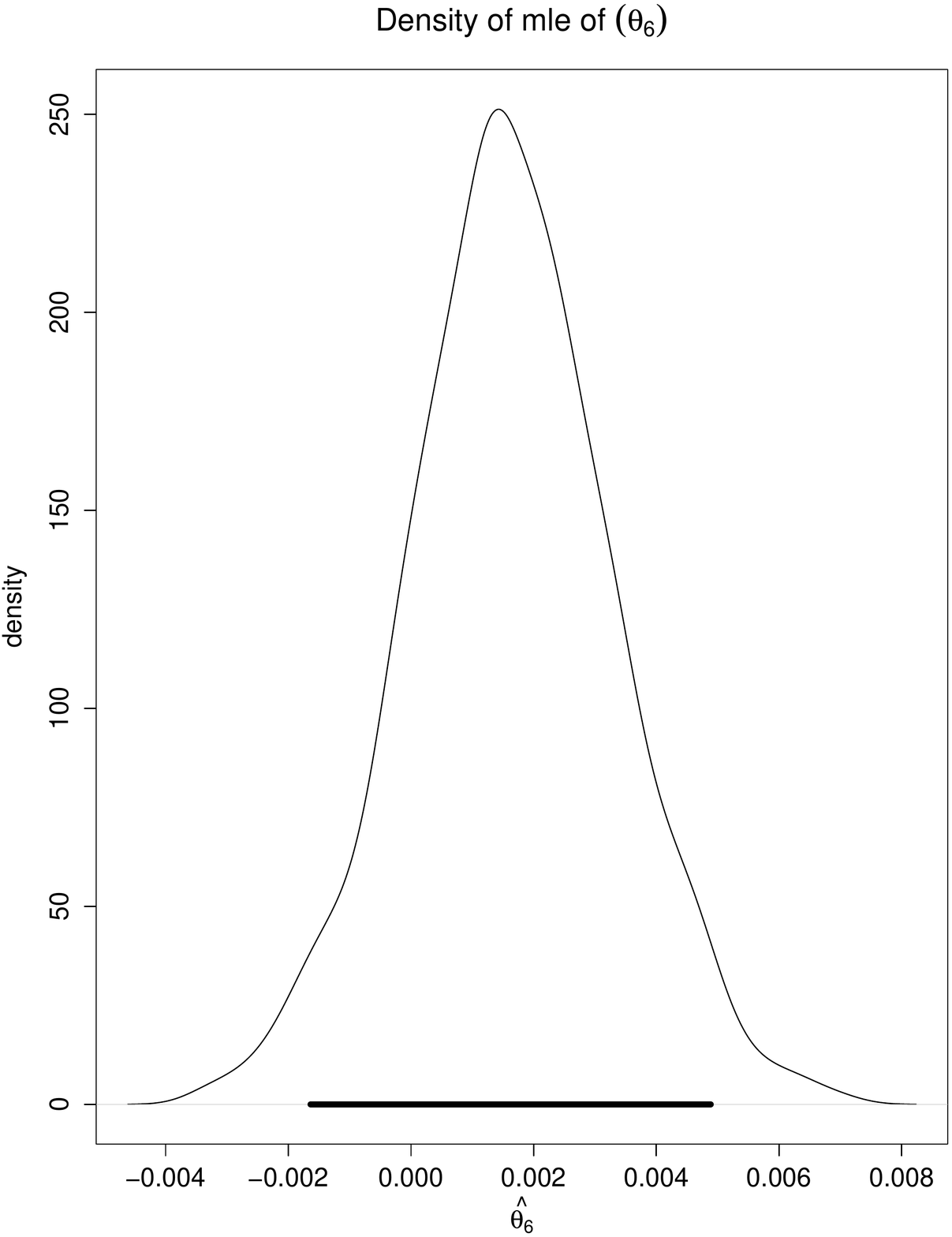}
\label{fig:sfig_r6}
\end{subfigure}
\label{fig:real_c_15}
\caption{Distribution of $MLE$s for the real data.}
\end{figure}

We first obtain the $MLE$s of the $6$ parameters $\theta_j$ for $j=1,\ldots,6$ by the block-relaxation algorithm described by Algorithm \ref{algo:block}, in Section \ref{subsec:C_N20}. 
In this real data set up, the process starts with the initial value $\theta_j=1$ for $j=1,\ldots,6$ (our experiments with several other choices demonstrated practicability of those
choices as well) and in step (3) of Algorithm \ref{algo:block}, the distance is taken as $0.1$ instead of $10^{-5}$. 
Then taking the $MLE$s as the value of the parameters $\theta_j$ for $j=1,\ldots,6$ we perform the parametric bootstrap method where we generate the data 1000 times and with 
respect to each data set, obtain the $MLE$s of the six parameters by the block-relaxation method as already mentioned.
Figure \ref{fig:real_c_15} shows the distribution of $MLE$s (denoted by $\hat\theta_j$ for $j=1,\ldots,6$)  where the respective $95\%$ confidence intervals are shown in bold lines.
Among the covariates, $c_3$, that is, the US dollar exchange rate, seems to be less significant compared to the others, since the distribution of the $MLE$ of the associated coefficient, $\theta_3$,
has highest density around zero, with small variability, compared to the other coefficients. Also note that the distribution of $\hat\theta_6$ is highly concentrated around zero,
signifying that the $X_i(t)$ term in the drift function of (\ref{eq:sde4_appl}) is probably redundant.

\subsection{Posterior Distribution of the parameters}
\label{subsec:B_R15}

In the Bayesian approach all the set up regarding the real data is exactly the same as in Section \ref{sec:truedata}, that is, each data is driven by the $SDE$s (\ref{eq:sde4_appl}) 
where the covariates $c_j$ for $j=1,\ldots,3$ are already mentioned in that section.
In this case, we consider the priors for the $6$ parameters to be independent normal with mean zero and variance $100$.
Since in real data situations the parameters are associated with greater uncertainties compared to simulation studies, somewhat vague prior as we have chosen here, as opposed
to that in the simulation study case, makes sense.

\begin{figure}
\begin{subfigure}
\centering
\includegraphics[height=6cm,width=5cm]{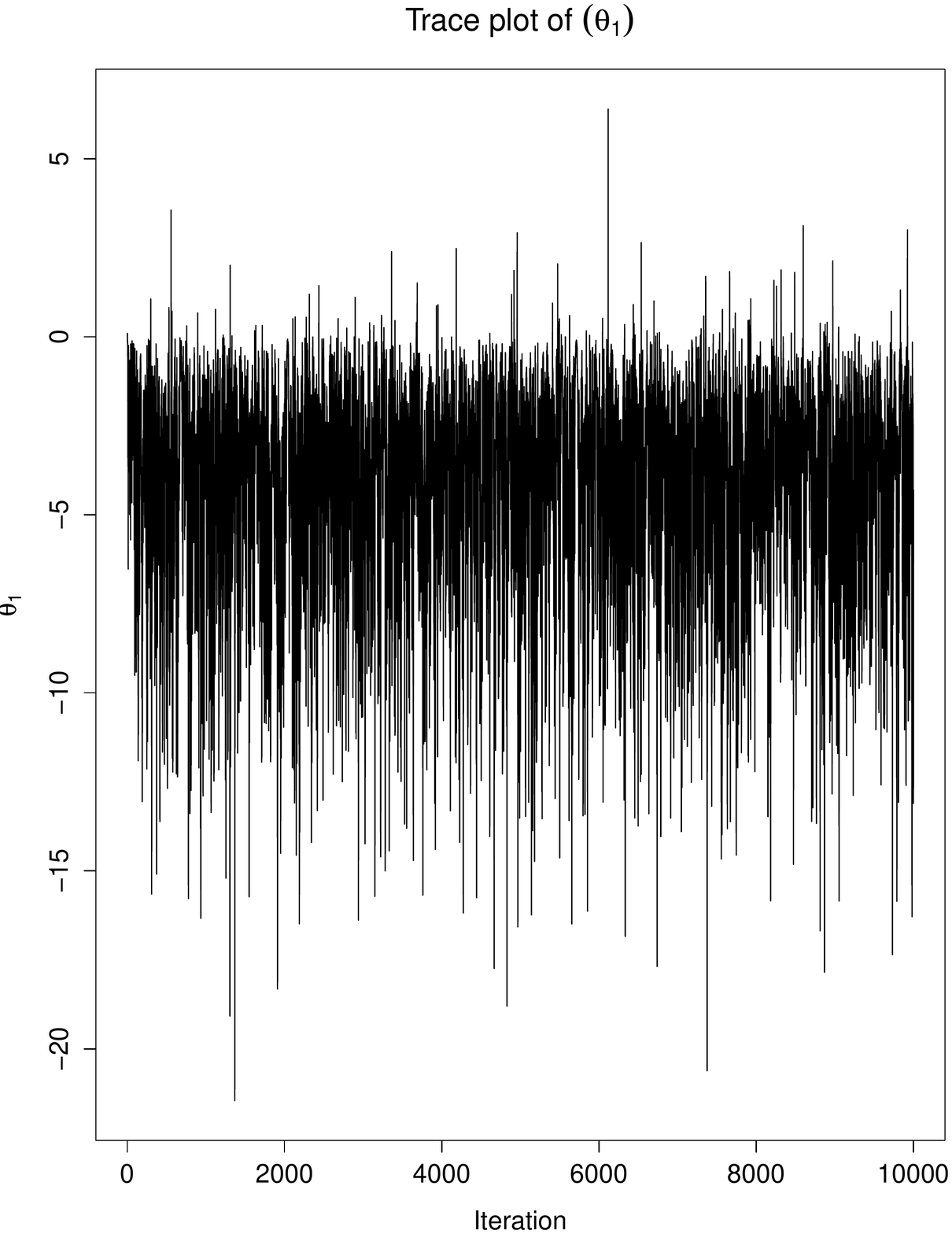}
\label{fig:sfig_t1}
\end{subfigure}
\begin{subfigure}
\centering
\includegraphics[height=6cm,width=5cm]{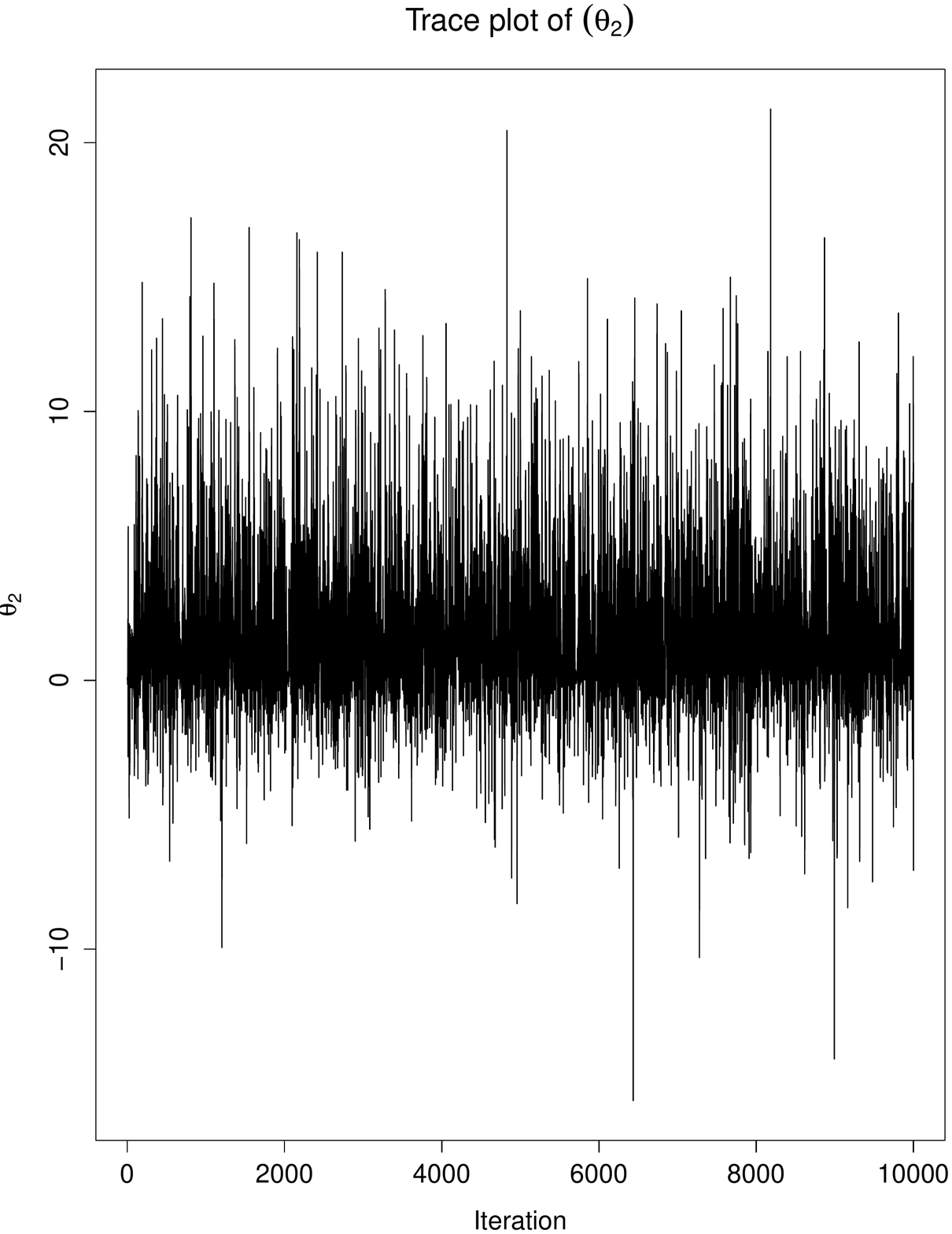}
\label{fig:sfig_t2}
\end{subfigure}
\begin{subfigure}
\centering
\includegraphics[height=6cm,width=5cm]{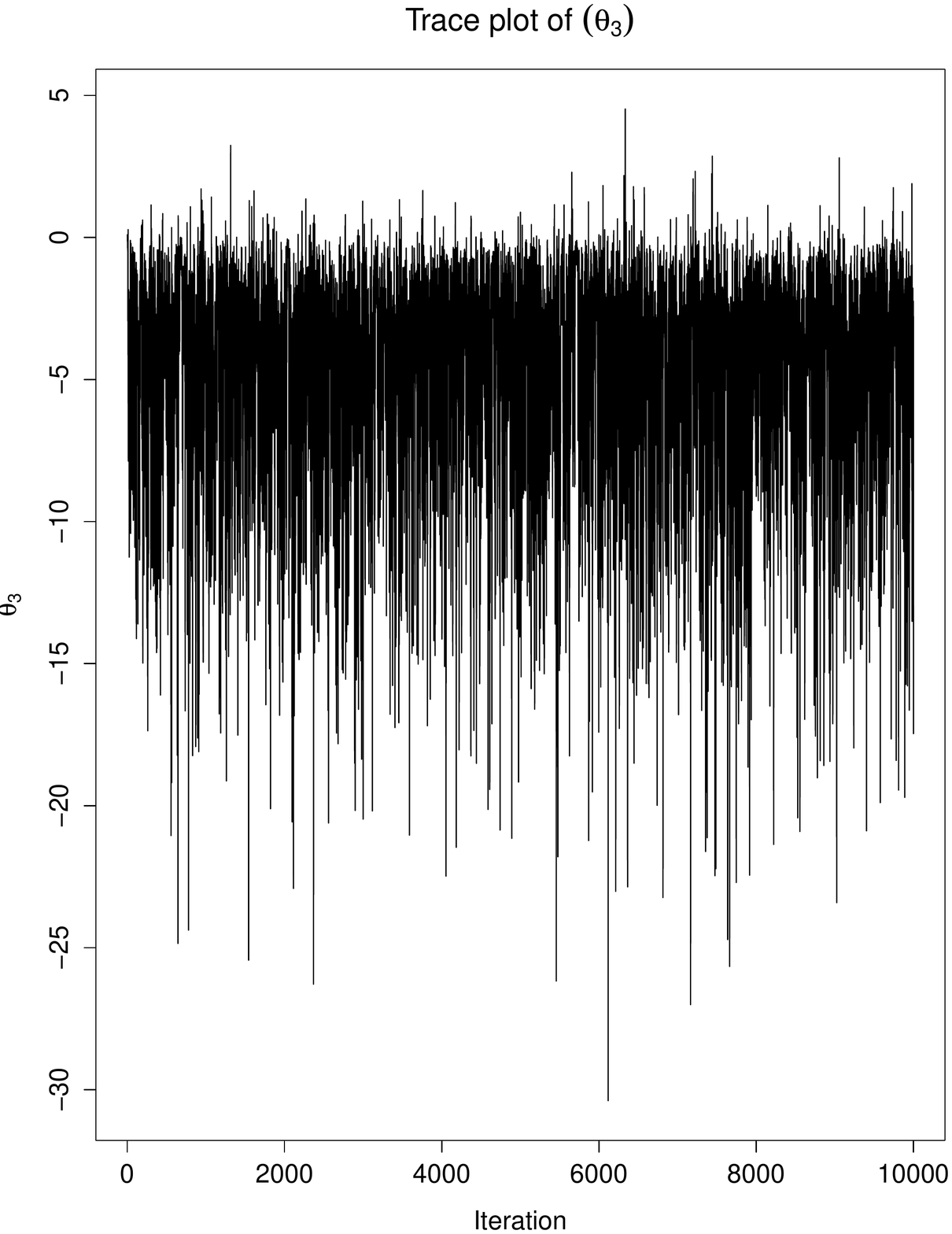}
\label{fig:sfig_t3}
\end{subfigure}
\begin{subfigure}
\centering
\includegraphics[height=6cm,width=5cm]{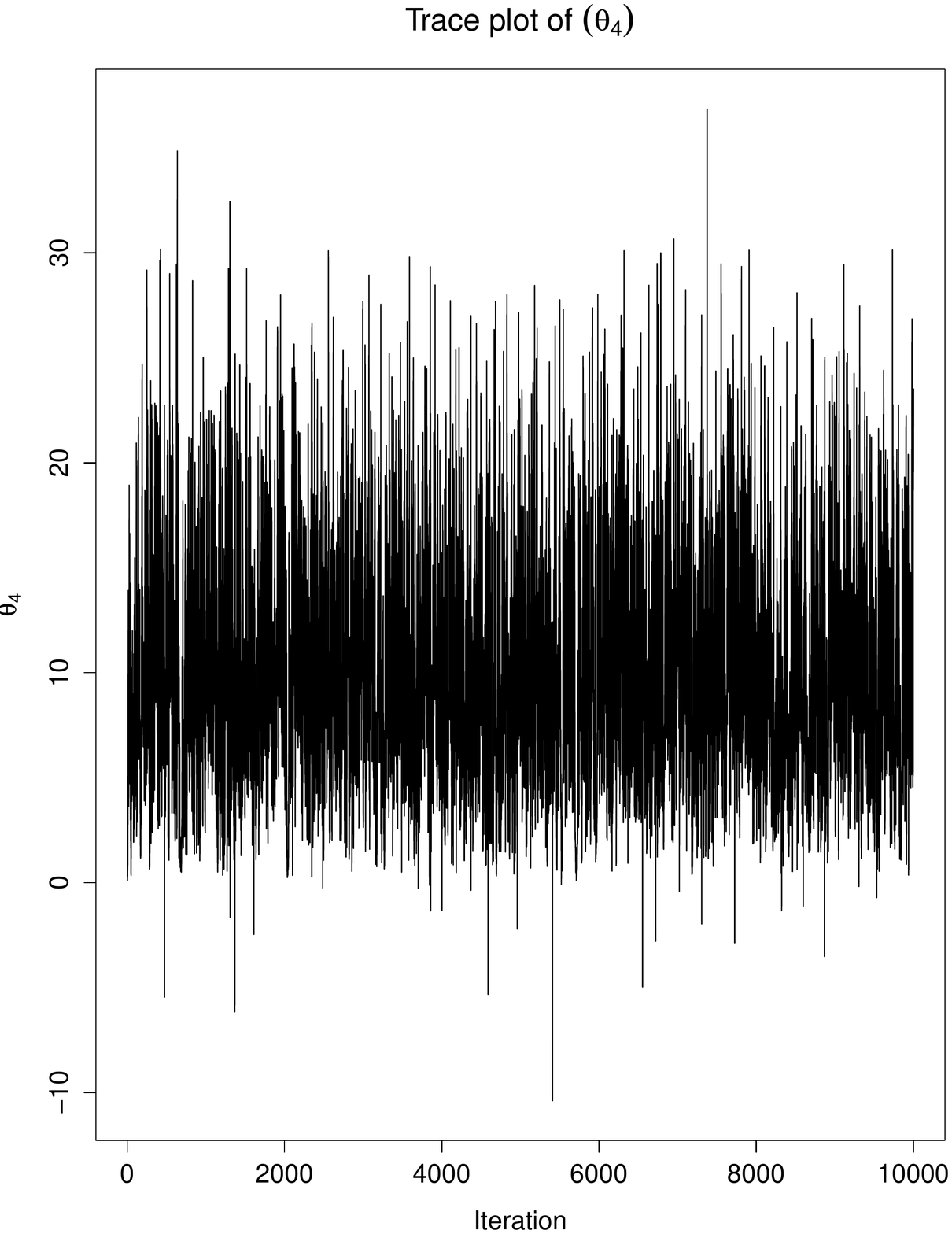}
\label{fig:sfig_t4}
\end{subfigure}
\begin{subfigure}
\centering
\includegraphics[height=6cm,width=5cm]{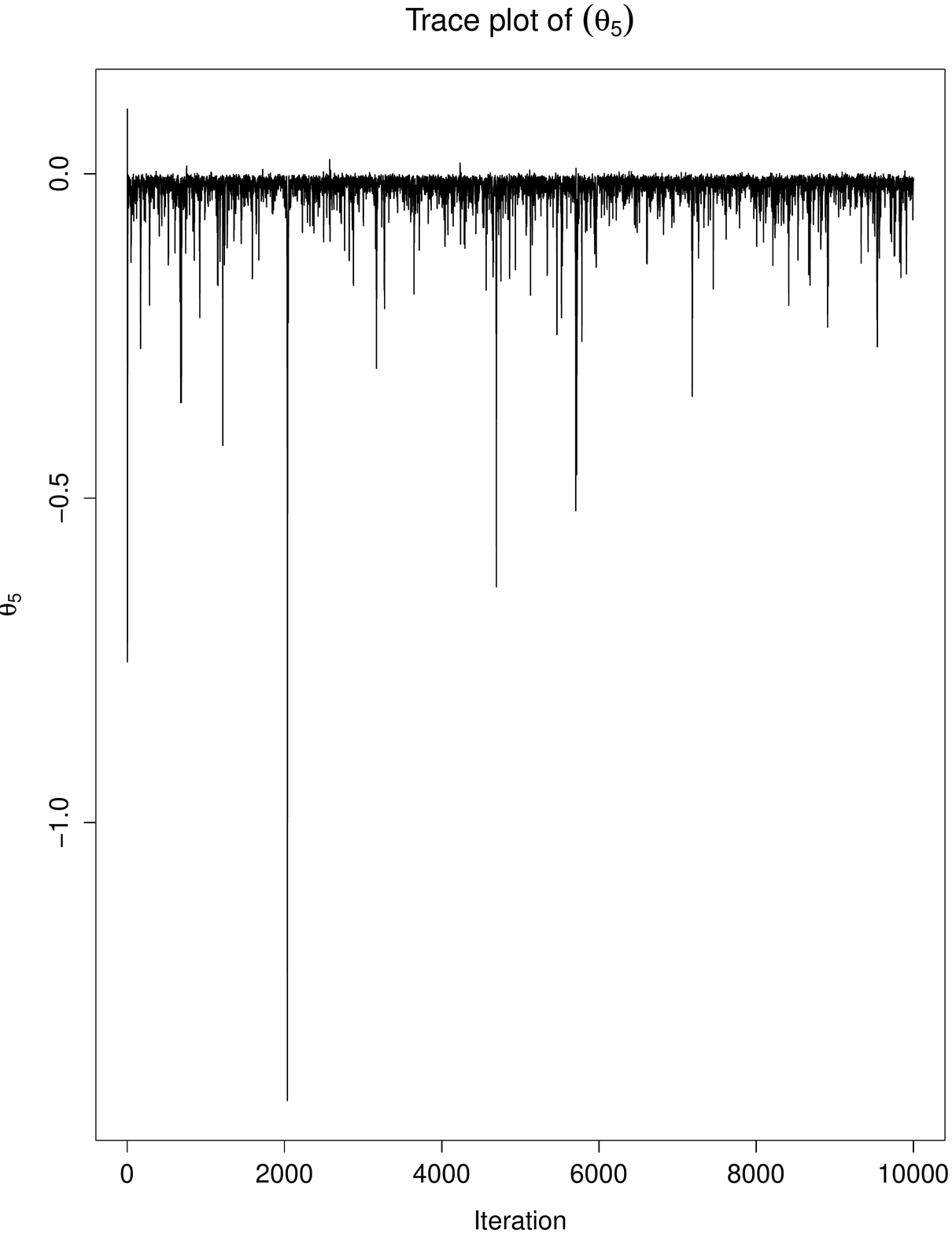}
\label{fig:sfig_t5}
\end{subfigure}
\begin{subfigure}
\centering
\includegraphics[height=6cm,width=5cm]{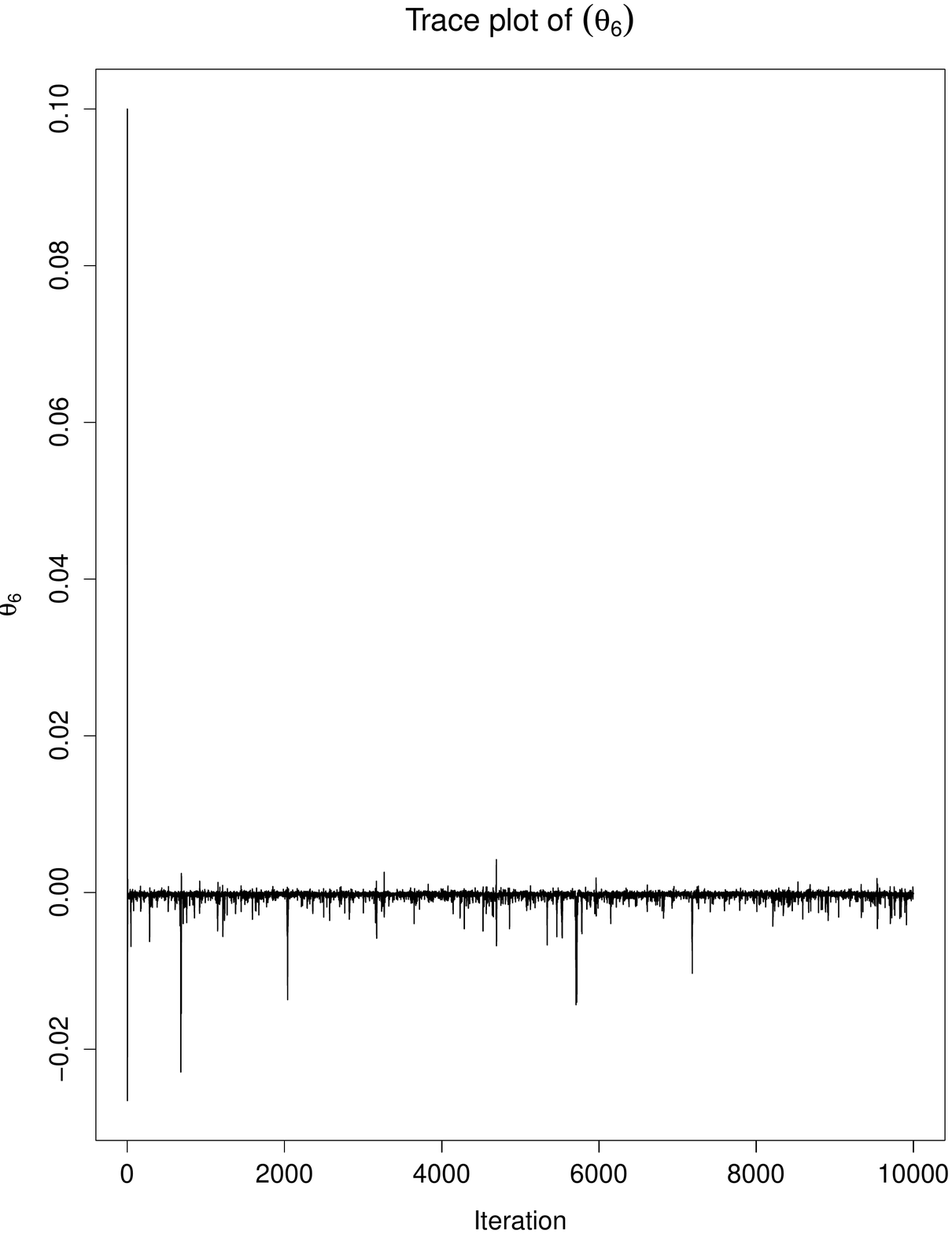}
\label{fig:sfig_t6}
\end{subfigure}
\label{fig:real_trace_15}
\caption{Trace Plot of the Parameters}
\end{figure}

The greater uncertainty in the parameters in this real data scenario makes room for more movement, and hence, better mixing of MCMC samplers such as Gibbs sampling, 
in contrast with that in simulation studies. As such, our application of Gibbs sampling, were the full conditionals are normal distributions with appropriate means and variances,
yielded excellent mixing. Although we chose the initial values as $\theta_j=0.1$; $j=1,\ldots,6$, other choices also turned out to be very much viable.
We perform $100000$ Gibbs sampling iterations to obtain our required posterior distributions of the $6$ parameters.
Figure \ref{fig:real_trace_15} shows the trace plots of the $6$ parameters associated with $10000$ thinned samples obtained by plotting the output of every $10$-th iteration. 
We emphasize that although we show the trace plots of only $10000$ Gibbs sampling realizations to reduce the file size, our inference is based on all the $100000$ realizations.
From the trace plots, convergence of the posterior distributions of the parameters is clearly observed.

\begin{figure}
\begin{subfigure}
\centering
\includegraphics[height=6cm,width=5cm]{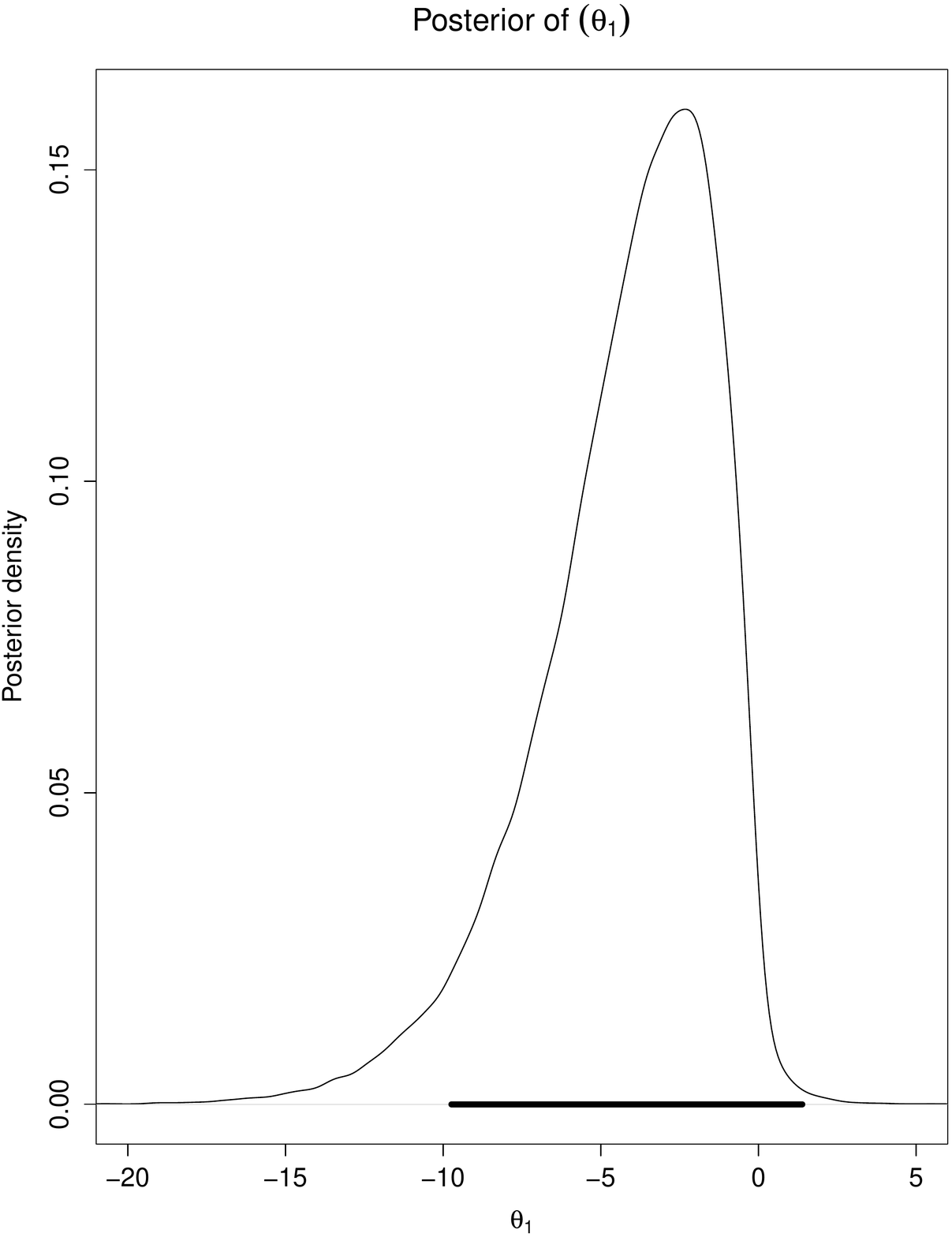}
\label{fig:sfig_r7}
\end{subfigure}
\begin{subfigure}
\centering
\includegraphics[height=6cm,width=5cm]{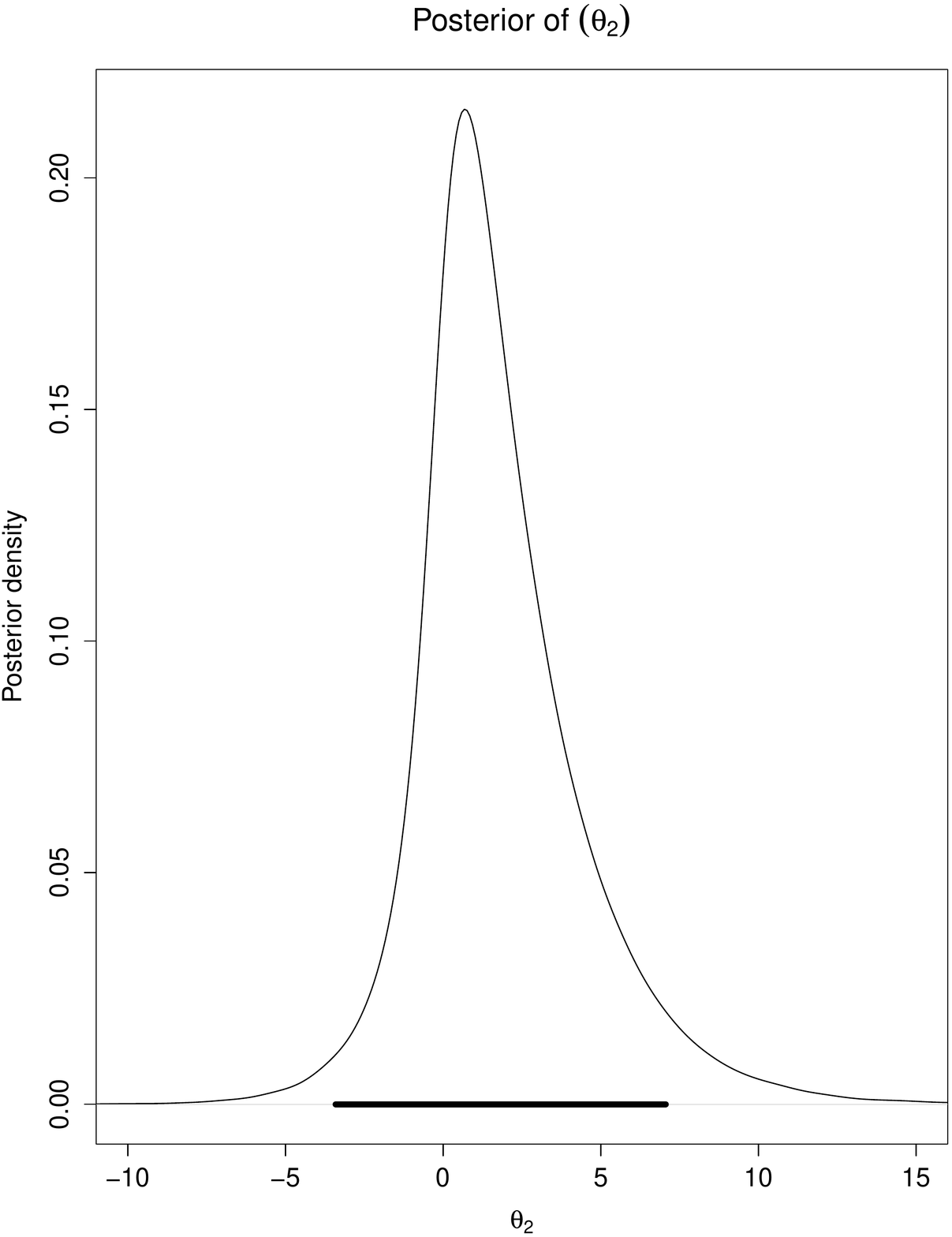}
\label{fig:sfig_r8}
\end{subfigure}
\begin{subfigure}
\centering
\includegraphics[height=6cm,width=5cm]{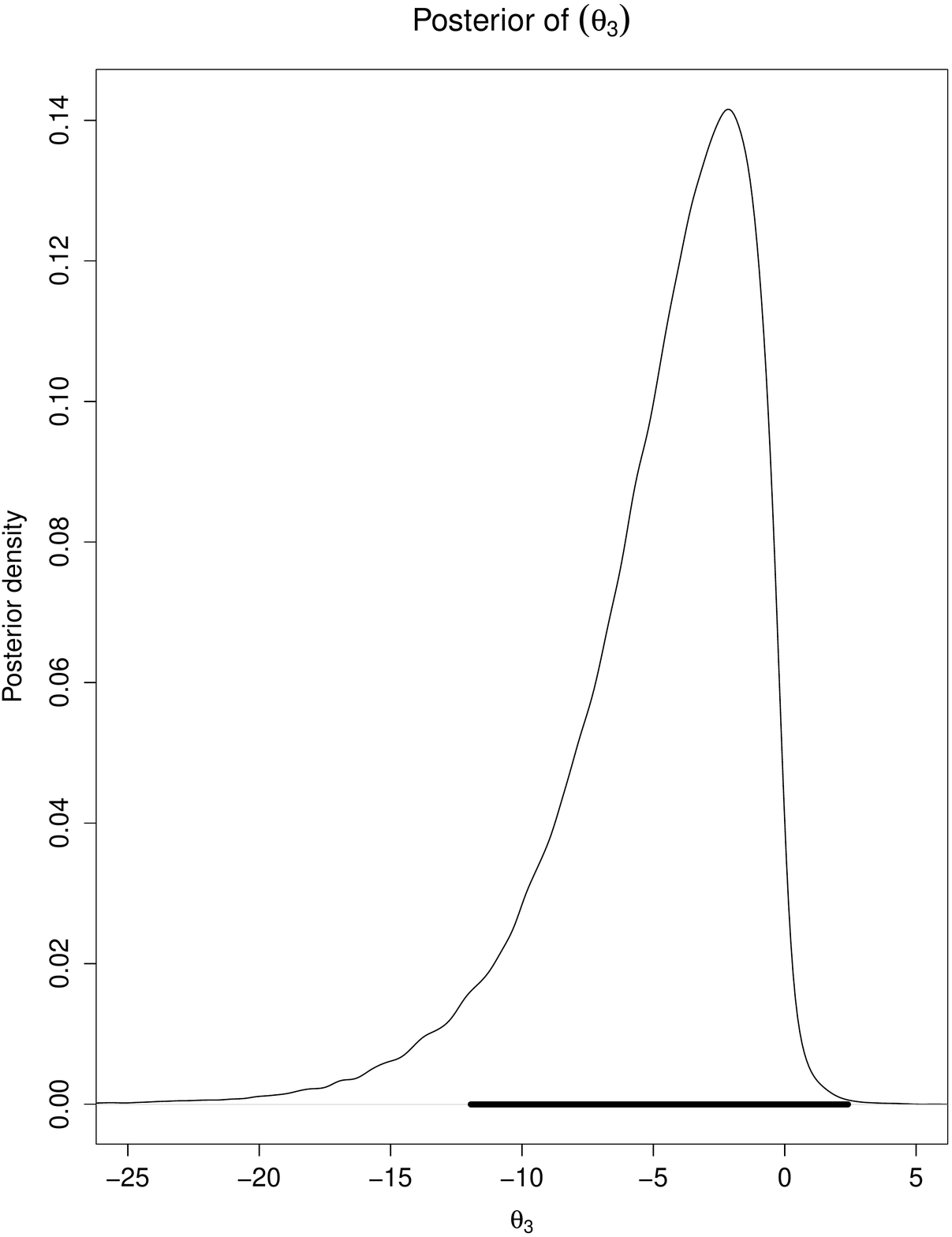}
\label{fig:sfig_r9}
\end{subfigure}
\begin{subfigure}
\centering
\includegraphics[height=6cm,width=5cm]{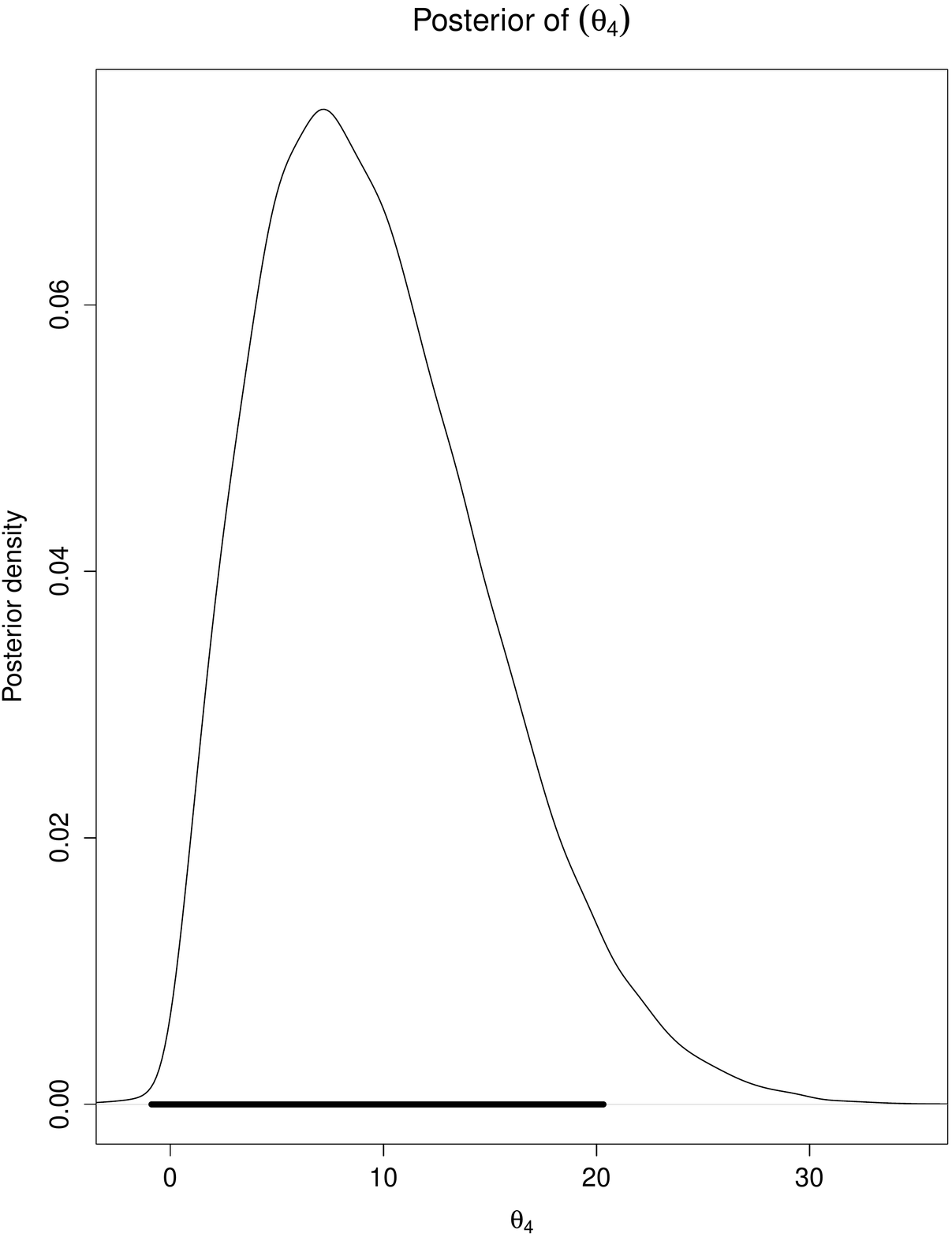}
\label{fig:sfig_r10}
\end{subfigure}
\begin{subfigure}
\centering
\includegraphics[height=6cm,width=5cm]{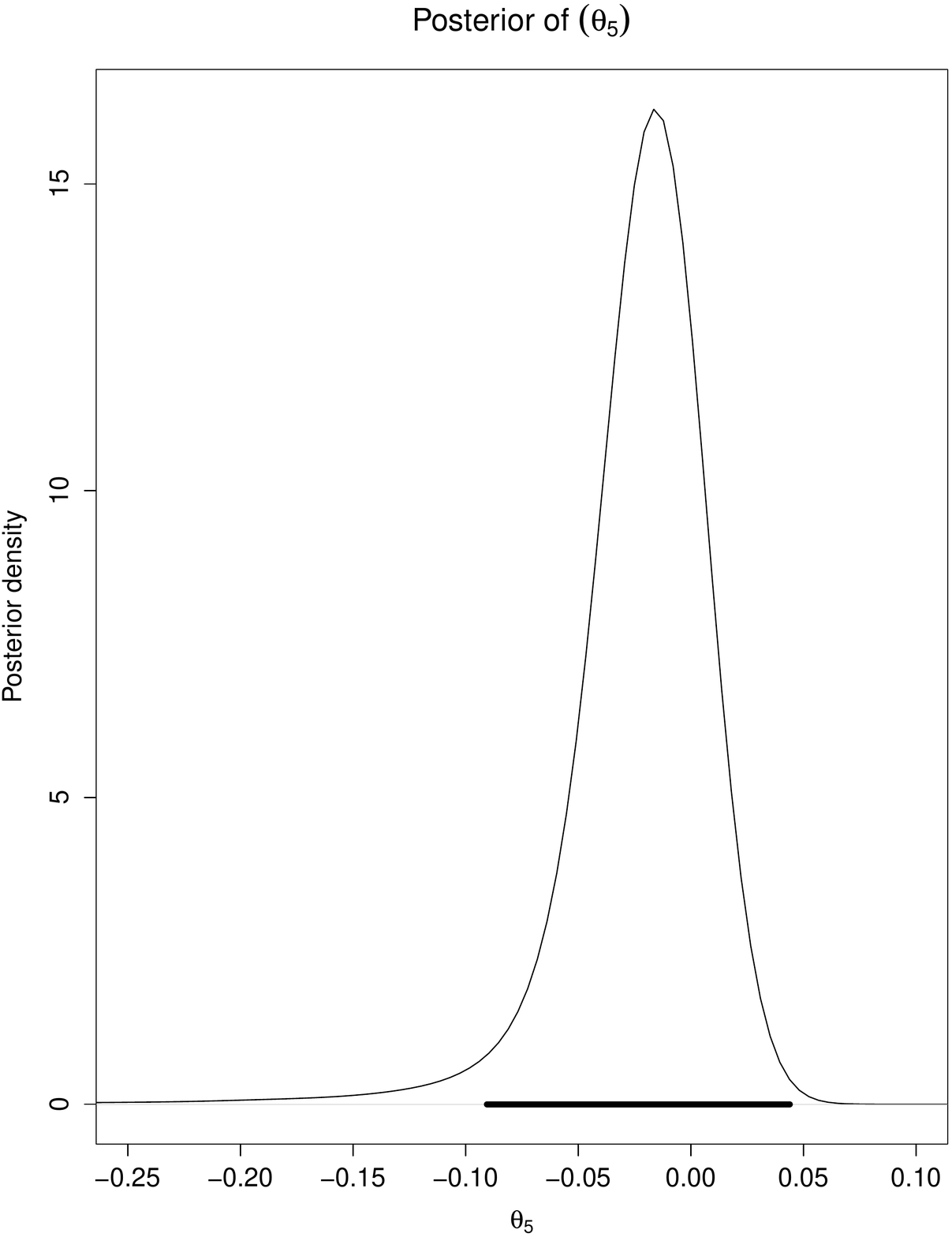}
\label{fig:sfig_r11}
\end{subfigure}
\begin{subfigure}
\centering
\includegraphics[height=6cm,width=5cm]{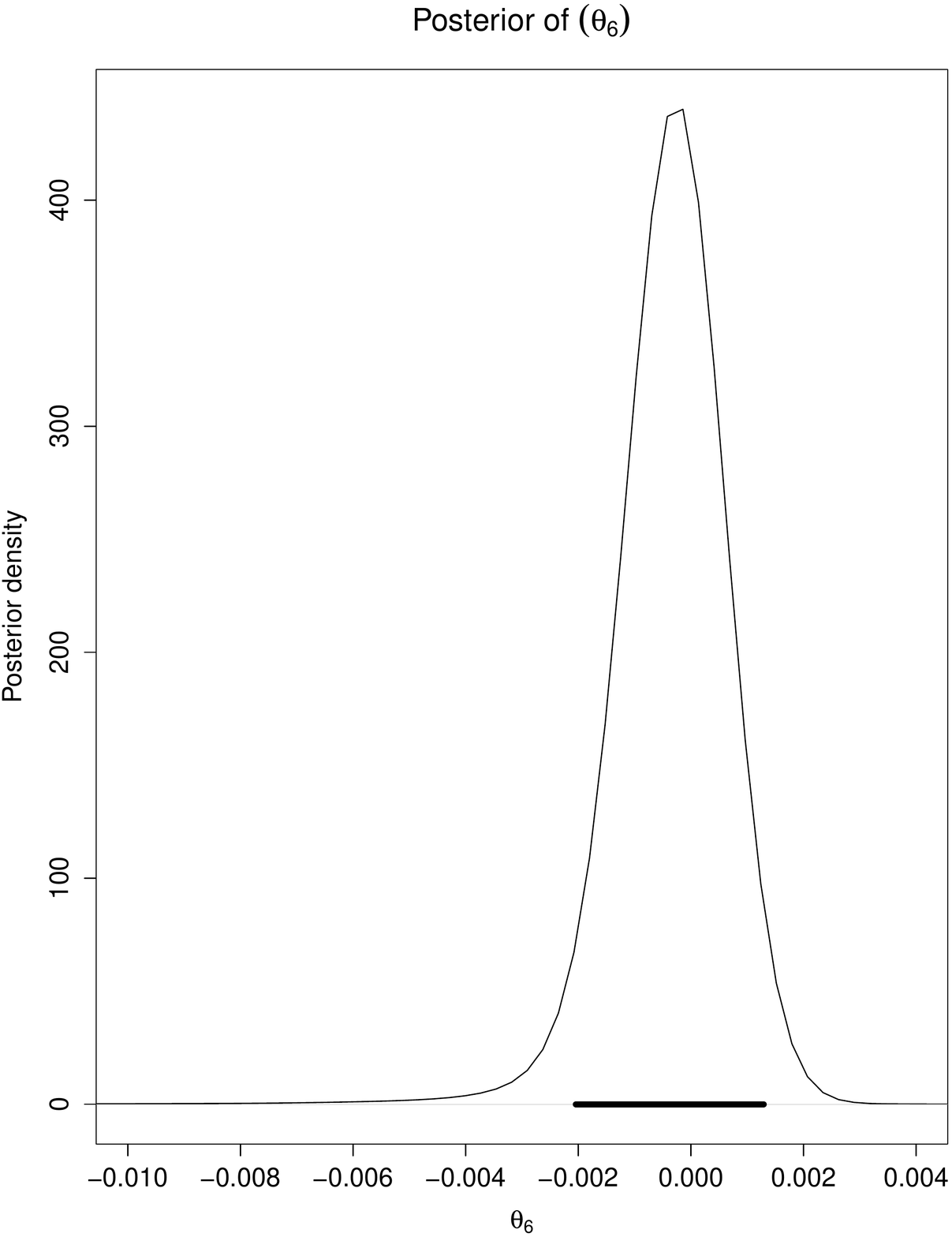}
\label{fig:sfig_12}
\end{subfigure}
\label{fig:real_b_15}
\caption{Posterior Distributions of the Parameters for real data}
\end{figure}

Figure \ref{fig:real_b_15} displays the posterior densities of the 6 parameters, where the $95\%$ credible intervals are indicated by bold lines.
The posterior distribution of $\theta_3$ is seen to include zero in the highest density region; however, unlike the distribution of the $MLE$ $\hat\theta_3$,
the posterior of $\theta_3$ has a long left tail, so that insignificance of $c_3$ is not very evident. The posterior of $\theta_6$ is highly concentrated around zero,
agreeing with the $MLE$ of $\theta_6$ that the term $X_i(t)$ in the drift function is perhaps redundant. Note that the posterior of $\theta_5$ also inclues zero
in its high-density region, however, it has a long left tail, so that the significance of $\theta_5$, and hence, of the overall drift function, is not ruled out.

\section{Summary and conclusion}
\label{sec:conclusion}

In $SDE$ based random effects model framework, \ctn{Maud12} considered 
the linearity assumption in the drift function given by $b(x,\phi_i) = \phi_ib(x)$, assuming $\phi_i$ to 
be Gaussian random variables with mean $\mu$ and variance $\omega^2$, and obtained a closed form 
expression of the likelihood of the above parameters. Assuming the $iid$ set-up, they proved convergence in 
probability and asymptotic normality of the maximum likelihood estimator of the parameters. 

\ctn{Maitra15} and \ctn{Maitra15b} extended their model by incorporating time-varying covariates in $\phi_i$
and allowing $b(x)$ to depend upon unknown parameters, but rather than inference regarding the parameters,
they developed asymptotic model selection theory based on Bayes factors for their purposes. 
In this paper, we developed asymptotic theories for parametric inference
for both classical and Bayesian paradigms under the fixed effects set-up, and provided relevant discussion
of asymptotic inference on the parameters in the random effects set-up. 

As our previous investigations (\ctn{Maitra14a}, \ctn{Maitra14b}, for instance), 
in this work as well we distinguished the non-$iid$ set-up from the $iid$ case, the latter corresponding
to the system of $SDE$s with same initial values, time domain, but with no covariates. However, as 
already noted, this still provides a generalization to the $iid$ set-up of \ctn{Maud12} through
generalization of $b(x)$ to $b_{\bbeta}(x)$; $\bbeta$ being a set of unknown parameters. Under suitable
assumptions we obtained strong consistency and asymptotic normality of the $MLE$ under the $iid$ set-up
and weak consistency and asymptotic normality under the non-$iid$ situation.
Besides, we extended our classical asymptotic theory to the Bayesian framework, for both
$iid$ and non-$iid$ situations. Specifically, we proved posterior consistency and asymptotic posterior
normality, for both $iid$ and non-$iid$ set-ups. 

In our knowledge, ours is the first-time effort regarding asymptotic inference, either classical or Bayesian, 
in systems of $SDE$s under the presence of time-varying covariates. 
Our simulation studies and real data applications, with respect to both classical and Bayesian paradigms, have revealed
very encouraging results, demonstrating the importance of our developments even in practical, finite-sample situations.



\section*{Acknowledgment}
We are sincerely grateful to the EIC, the AE and the referee whose constructive comments have
led to significant improvement of the quality and presentation of our manuscript.
The first author also gratefully acknowledges her CSIR Fellowship, Govt. of India.

\normalsize
\bibliographystyle{natbib}
\bibliography{irmcmc}

\end{document}